\documentclass[preprint,12pt]{elsarticle}

\usepackage[right=2.5cm,left=2.5cm,bottom=2.5cm,top=2.5cm]{geometry}
\usepackage[reqno]{amsmath}
\usepackage[utf8]{inputenc}
\usepackage{amssymb}
\usepackage{epsfig}
\usepackage{longtable}
\usepackage{xspace}
\usepackage{rotating}
\usepackage{appendix}
\usepackage{graphicx}
\usepackage{epstopdf}
\usepackage{color}
\usepackage{subfigure}
\usepackage{booktabs}
\usepackage{url}
\usepackage{times}
\usepackage{paralist}
\usepackage{cleveref}
\usepackage{parskip}
\usepackage{mathrsfs}
\usepackage{amsmath}
\usepackage{amsthm}
\usepackage{tikz}
\usepackage{caption}
\usepackage{graphicx}
\usepackage[algoruled,linesnumbered]{algorithm2e}
\usepackage{microtype}
\usepackage{subfigure}
\usepackage{graphicx}
\usepackage{rotating}
\usepackage{algorithmicx}
\usepackage{pgfplots}
\usepackage{pgfplotstable}

\newtheorem{property}{Property}

\newcommand{\argmax}{\operatornamewithlimits{argmax}}

\usetikzlibrary{arrows}
\usetikzlibrary{calc}
\usetikzlibrary{fit}                      
\usetikzlibrary{backgrounds}              
\usetikzlibrary{shapes.geometric}
\usetikzlibrary{shapes}
\usetikzlibrary{shapes.multipart}
\usetikzlibrary{decorations.pathreplacing}
\usetikzlibrary{graphs}
\usetikzlibrary{positioning}
\tikzset{main node/.style={circle,draw,minimum size=0.55cm,inner sep=0pt},}
\tikzset{dsatur node/.style={circle,draw,minimum size=0.6cm,inner sep=0pt},}
\tikzset{reduced node/.style={circle,draw,minimum size=1cm,inner sep=0pt},}
\tikzset{little node/.style={circle,draw,minimum size=0.1cm,inner sep=0pt},}
\tikzset{small node/.style={circle,fill=black,draw,minimum size=0.15cm,inner sep=0pt},}

\SetAlgoSkip{smallskip}
\DontPrintSemicolon
\crefname{algocf}{algorithm}{algorithms}
\Crefname{algocf}{Algorithm}{Algorithms}
\crefname{algocfline}{line}{lines}
\Crefname{algocfline}{Line}{Lines}
\crefname{AlgoLine}{line}{lines}
\Crefname{AlgoLine}{Line}{Lines}

\newlength{\mylengthleft}
\newlength{\mylength}
\usepackage{tikz}

\usetikzlibrary{arrows,shapes,positioning}
\usetikzlibrary{decorations.markings}
\tikzstyle arrowstyle=[scale=1.8]
\tikzstyle directed=[postaction={decorate,decoration={markings,
    mark=at position .65 with {\arrow[arrowstyle]{stealth}}}}]
\tikzstyle reverse directed=[postaction={decorate,decoration={markings,
    mark=at position .65 with {\arrowreversed[arrowstyle]{stealth};}}}]

    \tikzset{
mainLine/.style={thick},
secondLine/.style={black},
}

\def\BPP{{\rm BPP}}

\def\TBPP{{\rm TBPP}}
\def\TKP{{\rm TKP}}
\def\VPP{{\rm VPP}}
\def\SSP{{\rm SSP}}
\def\2DBPP{{\rm 2D-BPP}}
\def\CSP{{\rm CSP}}

\def\LP{{\rm LP}}
\def\ILPC{\text{\rm ILP$^c$}}
\def\ILPE{\text{\rm ILP$^e$}}
\def\LPE{\text{\rm LP$^e$}}
\def\LPC{\text{\rm LP$^c$}}

\journal{{}}

\begin{document}

\begin{frontmatter}
  \title{\Large {A Branch-and-Price Algorithm for the  {Temporal} Bin {Packing} Problem}}

  \author{Mauro Dell'Amico}
  \address{DISMI, University of Modena and Reggio Emilia,
Via Amendola 2, 42122 Reggio Emilia, Italy\\ mauro.dellamico@unimore.it}
  \author{Fabio Furini \footnote{Corresponding author}}
  \address{Universit\'e Paris Dauphine, PSL Research University, LAMSADE,  75016 Paris,  France\\ fabio.furini@dauphine.fr}
  \author{Manuel Iori }
  \address{DISMI, University of Modena and Reggio Emilia,
Via Amendola 2, 42122 Reggio Emilia, Italy \\ manuel.iori@unimore.it}

  \begin{abstract}
We study an extension of the classical Bin Packing Problem, where each item consumes the bin capacity during a given time window {that depends} on the item itself.  The problem asks for finding the minimum number of bins to pack all the items while respecting the bin capacity  at any time instant. {A polynomial-size formulation, an exponential-size formulation, and a number of lower and upper bounds} are studied. A branch-and-price algorithm for solving the exponential-size formulation is introduced. An overall algorithm combining the different methods is then proposed and tested trough extensive computational experiments.
\end{abstract}

\end{frontmatter}


\section{Introduction}\label{sec:intro}

The \textit{Bin Packing Problem} (\BPP) is one of the classical problems in combinatorial optimization and has been extensively studied in the literature, see, e.g., \citet{DIM16}. Given a large number of identical bins of capacity ${W \in \mathbb{Z}_{+}}$ and a set $N=\{1, 2, \ldots, n\}$ of items, where each item $j\in N$ is associated with an integer weight ${w_j}\leq W$, the \BPP\ asks to {pack} all the items {into} the minimum number of bins {without exceeding the capacity}.

In this paper, we study a natural generalization of the \BPP\ called  {the} \textit{Temporal Bin Packing Problem} (\TBPP). In the \TBPP, a feasible assignment of the items to the bins must be computed over a \emph{{discretized} time horizon} ${\widehat{T}=\{0,1,2,\dots,|\widehat{T}|\}}$, {of total length $|\widehat{T}|$}. If a bin is selected, its capacity is a renewable resource that is available at any time unit in the horizon. Each item $j\in N$  consumes  the bin capacity during a {given} time window $[s_j, t_j)$, with $0 \leq s_j < t_j < {|\widehat{T}|}$. The {integer} input parameters {$s_j \in \mathbb{Z}_{+}$} and {$t_j \in \mathbb{Z}_{+}$} represent the \emph{starting time} and the \emph{ending time} of an item, respectively. {As for the \BPP, in the \TBPP\  each item must be assigned to a unique bin where it remains for its entire time window.}
{The \TBPP\ asks to  pack all the items  into the minimum number of bins
so that the bin capacity is never exceeded at any unit of time.}
The problem is \emph{strongly} NP-hard, {because} its restriction obtained by setting $|\widehat{T}|=1$ {boils down to a} \BPP, which is well known to be strongly NP-hard.

The \TBPP\  adds  a temporal dimension to the classical \BPP, thus making the problem  very challenging to solve in practice. A related difficult problem is the \emph{Vector Packing Problem} (\VPP), see, e.g., \citet{HGI18}, where each bin has $k$ capacities $W_1,\dots,W_k$ and each item $j$ is associated with a vector of $k$ weights $w_j=(w_{j1}, \dots, w_{jk})$. A feasible \VPP\ solution consists in packing all items in the bins so that the capacity is respected for all the $k$ dimensions. We will show in Section \ref{sec:Models} that the \TBPP\ is a special case of the \VPP.

Other related problems are: {i)} the \emph{Two-Dimensional} BPP (\2DBPP),  see, e.g. \citet{PR07}, in which both items and bins are rectangles and the aim is to pack all items without overlapping in the minimum number of bins; {and ii)} the \emph{Temporal Knapsack Problem} (TKP), see, e.g., \citet{CFM13},  in which items {also have} a profit and the aim is to find a subset of items of maximum profit that fits into a single bin. These problems, as well as other interesting related problems, are discussed in more detail in Section \ref{sec:literature}.

The \TBPP\ finds applications in many fields, including logistics, {healthcare}, production and warehouse management. Consider, for example, the production field: each item can be interpreted as a task ({or a product}) that must attain a given production rate (equal to its weight) in each {time unit} of a given time window; each bin can be seen as a production plant that can be used for the allocation of tasks; minimizing the number of bins consequently implies minimizing the number of production plants that are used for the tasks (see, e.g., \citet{Angelelli14} for a related problem).

The remainder of the paper is organized as follows.
Section \ref{sec:literature} describes the related literature. Section \ref{sec:Models} presents two  \TBPP\  mathematical models, the former having polynomial size and the latter exponential. Sections \ref{sec:lower} and \ref{sec:ub} provide, respectively, lower and upper bounds for the problem, whereas Section \ref{sec:BP} describes a branch-and-price algorithm that solves the exponential-size model.  All algorithms and models are computationally tested in Section \ref{sec:computational} and conclusions are drawn in Section \ref{sec:conclusions}.
\section{Literature and related problems}\label{sec:literature}
The \BPP\ is one of the most widely studied problems in the combinatorial optimization field. A number of surveys and annotated bibliographies have been consequently proposed during the years to describe the main techniques that have been developed for its solution. Such techniques are either focused on the \BPP\ or on its reformulation known as the \emph{Cutting Stock Problem} (\CSP), where all items having same weight are grouped together into {item types}.

Useful classifications have been provided by \citet{WHS07}, who presented a typology of cutting and packing problems based on  detailed categorization criteria, and \citet{CC07b}, who introduced a four-field classification scheme aimed at highlighting the
main theoretical results in the area. A few years later, \citet{CCGMV13} presented an overview of approximation algorithms for the BPP and a number of its variants, and classified all references according to \citet{CC07b}. \citet{V02} presented a survey with a focus on the most popular {\em Linear Programming} (LP) methods for the \BPP\ and the \CSP.

Recently, \citet{DIM16} reviewed the most important mathematical models and algorithms developed for the exact solution of the \BPP\ and the \CSP, and experimentally evaluated the performance of the main available software tools. The extensive results obtained, together with the input benchmark instances addressed, have been gathered together and made available on-line at the \emph{Bin Packing Problem Library}, as shown in \citet{DIM18}. {Exact algorithms that appeared after \citet{DIM18} are the iterative aggregation and disaggregation method by \citet{CHMVA17}, the improved \emph{reflect} formulation of \citet{DI18}, and the branch-and-price by \citet{Wei2018}.}

A number of problem extensions have been proposed during the years. We believe it is worth describing the main results that have been obtained on those variants that are close to the \TBPP.

The previously mentioned VPP has been the object of several interesting researches. \citet{CT01} focused on the case with $k=2$ dimensions, providing effective heuristics and a few exact algorithms, the most effective one based on column generation.
\citet{AVCR14} implemented several dual-feasible functions and fast lower bounding techniques.
\citet{Brandao2016} used pseudo-polynomial arc-flow models and managed to reduce their size through the use of graph reduction techniques. Very recently, \citet{HGI18} proposed efficient stabilized branch-and-price algorithms. Their column-generation sub-problem is a multidimensional knapsack problem (see, e.g., \citet{DDIM19}) either binary, bounded, or unbounded, that they solved as a shortest path problem with resource constraints.

Extensions in which items and bins are boxes in $d$ dimensions have also been intensively studied. Most of the works on these problems focused on the case where $d=2$, solving the \2DBPP. The aim of the \2DBPP\ is to pack all items into the minimum numbers of bins without overlapping. A review of some of the methods to solve the \BPP\ and the \2DBPP\ was given in the early nineties by \citet{HS91}. Later on, surveys on the \2DBPP\ and on some of its relevant variants were proposed by \citet{LMM02,LMV02,LMMV10}.
Recent relevant results on the \2DBPP\ have been obtained, among others, by \citet{PR07}, who developed an efficient exact algorithm based on column generation and constraint programming, and by \citet{MM18}, who proposed a list of lower bounding techniques.

The problem of \emph{interval scheduling with a resource constraint} (ISRC) was presented by \citet{Angelelli11}. The ISRC is a scheduling problem where jobs have to be processed by parallel identical machines. Similarly to the items in \TBPP, each job in the ISRC has fixed start and finish time, as well as a resource consumption (i.e., a weight).
\citet{Angelelli11} focused on the recognition version of the ISRC, and proved that deciding whether an instance has a feasible solution is strongly NP-complete even when the resource capacity of the machines is fixed to any value greater than or equal to two.
A few years later, \citet{Angelelli14} studied the optimization version of the ISRC, whose objective is a weighted function that depends on the assignment of jobs to machines. They proposed a column generation scheme, as well as greedy and restricted enumeration heuristics, {and extensively tested them on a number of instances.}

Another relevant problem is the \emph{BPP with Contiguity Constraints} (BPPC). In the BPPC, a certain number of copies might exist for an item, and all copies should be packed in consecutive bins.
Starting from \citet{MMV03}, the BPPC has been used as a relaxation for two-dimensional cutting and packing problems, either within branch-and-bound algorithms (see, e.g., \citet{APT09} and \citet{BR13}) or in combinatorial Benders decompositions (see, e.g., \citet{CDI14} and \citet{DIM17}).

It is well-known that the \BPP\ can be solved by a \emph{Dantzig-Wolfe}  reformulation in which each subproblem is a one-dimensional \emph{Knapsack Problem} (\citet{Gilmore1961,Gilmore1963}). The same result holds for the \TBPP, with the relevant difference that the subproblem is a TKP.
The TKP has received a fair amount of attention in the recent combinatorial optimization literature. The problem was formally introduced in
\citet{BFHMTU05} to model resource allocation problems in the context of sparse resources, such as communication bandwidth of computer memory. \citet{CFM13} were the first to solve the TKP with a Dantzig-Wolfe reformulation, using two variants of a branch-and-price algorithm in which subproblems are either associated with groups of capacity constraints or with single capacity constraints, and showing that the former variant performs much better than the latter.
\citet{GI17} provided improved computational results by producing stabilized column generation algorithms based on the use of  dual-optimal inequalities, i.e., inequalities that are fulfilled by at least one of the dual optimal solutions and can thus be used to reduce the search space (see \citet{V05} and \citet{ADV06} for further details on this type of techniques).
We also mention that \citet{CFMT16} solved the TKP by using a so-called recursive Dantzig-Wolfe reformulation, which uses the  reformulation not only for solving the original master problem, but also for recursively solving the pricing sub-problems.
\section{Mathematical models}\label{sec:Models}
In order to derive suitable models for \TBPP, we start by showing that the \TBPP\ can be modeled by considering only a polynomial number of time instants. Despite the capacity requirements being defined on the entire time horizon, a weight variation may arise only at the starting time of an item. Therefore, it is sufficient to satisfy the capacity restrictions at the $n$  starting times of the items. Given an item $j\in N$, let us define $\overline S_j:= \{\ell\in N : s_\ell \le s_j \mbox{ and } t_\ell > s_j\}$ as the set of active items at time $s_j$ (note that $j\in\overline S_j$).
As all items in $\overline S_j$ are active at the same time instant,  a capacity constraint must be imposed for these items.
Moreover, if $\overline S_j\subseteq \overline S_k$, then the associated capacity constraint at time $s_j$ is dominated by that of time $s_k$. Let us define $\overline T=\{t\in N:\overline S_t\nsubseteq \overline S_k, \forall k\in N\}$ as  the index set of all the \emph{non-dominated} constraints (or non-dominated sets). To model the problem, it is enough  to consider the capacity usage at each $t\in \overline T$.

To simplify the presentation we compact the indices in $\overline T$ by shifting them into the set $T=\{1,\dots,|\overline T|\}$ and we rename the corresponding non-dominated sets as $S_1,\dots,S_{|T|}$. We call $t\in T $ a \textit{time step}. W.l.o.g., let us also suppose that the items are sorted by non-decreasing starting times.

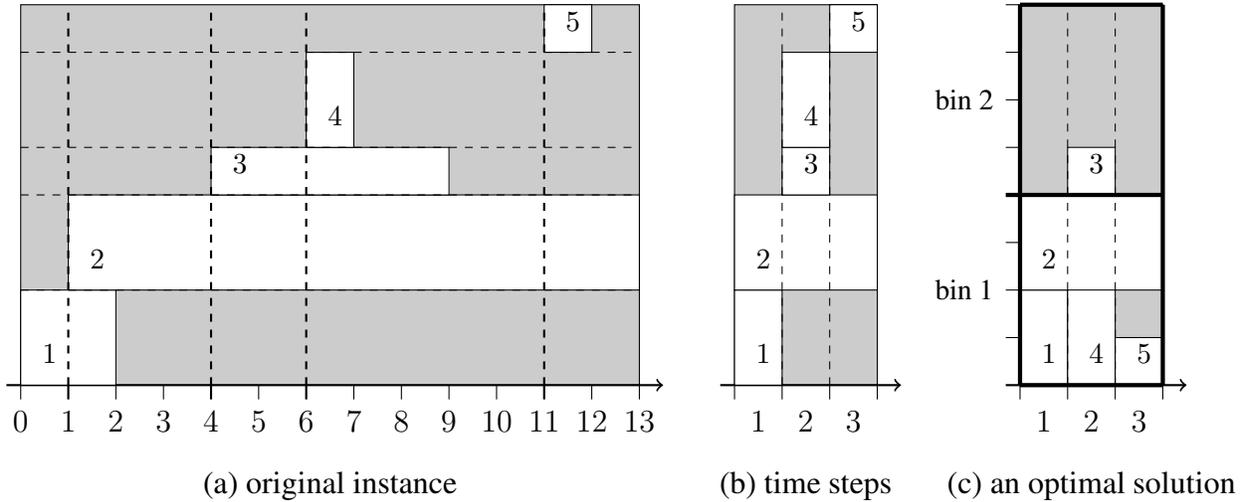
\begin{figure}[ht]
\begin{center}
\ifx\JPicScale\undefined\def\JPicScale{6.5}
\unitlength \JPicScale mm

\begin{tikzpicture}
\tikzset{x=18,y=18}


\draw[fill=black!20!white] (1,0) rectangle (14,8);

\draw[thick,->] (0.7,0)--(14.5,0)node[below right=0.2cm]{};
\draw (7.5,0) node[below =1.cm] {(a) original instance};

\draw[-] (1,0)--(1,-0.3);
\draw[thick,-] (2,0)--(2,-0.3);
\draw[-] (3,0)--(3,-0.3);
\draw[-] (4,0)--(4,-0.3);
\draw[thick,-] (5,0)--(5,-0.3);
\draw[-] (6,0)--(6,-0.3);
\draw[thick,-] (7,0)--(7,-0.3);
\draw[-] (8,0)--(8,-0.3);
\draw[-] (9,0)--(9,-0.3);
\draw[-] (10,0)--(10,-0.3);
\draw[-] (11,0)--(11,-0.3);
\draw[thick,-] (12,0)--(12,-0.3);
\draw[-] (13,0)--(13,-0.3);
\draw[-] (14,0)--(14,-0.3);

\draw (1,0) node[below =0.2cm] {$0$};
\draw (2,0) node[below =0.2cm] {$1$};
\draw (3,0) node[below =0.2cm] {$2$};
\draw (4,0) node[below =0.2cm] {$3$};
\draw (5,0) node[below =0.2cm] {$4$};
\draw (6,0) node[below =0.2cm] {$5$};
\draw (7,0) node[below =0.2cm] {$6$};
\draw (8,0) node[below =0.2cm] {$7$};
\draw (9,0) node[below =0.2cm] {$8$};
\draw (10,0) node[below =0.2cm] {$9$};
\draw (11,0) node[below =0.2cm] {$10$};
\draw (12,0) node[below =0.2cm] {$11$};
\draw (13,0) node[below =0.2cm] {$12$};
\draw (14,0) node[below =0.2cm] {$13$};



\draw[fill=white] (1,0) rectangle (3,2);
\draw (1,0) node[above right =0.2cm] {\small $1$};

\draw[fill=white] (2,2) rectangle (14,4);
\draw (2,2) node[above right =0.2cm] {\small  $2$};

\draw[fill=white] (5,4) rectangle (10,5);
\draw (5,4) node[above right =0.2cm] {\small $3$};

\draw[fill=white] (7,5) rectangle (8,7);
\draw (7,5) node[above right =0.2cm] {\small $4$};

\draw[fill=white] (12,7) rectangle (13,8);
\draw (12,7) node[above right =0.2cm] {\small $5$};

\draw[dashed] (1,2)--(14,2);
\draw[dashed] (1,4)--(14,4);
\draw[dashed] (1,5)--(14,5);
\draw[dashed] (1,7)--(14,7);

\draw[dashed, thick] (2,0)--(2,8);
\draw[dashed, thick] (5,0)--(5,8);
\draw[dashed, thick] (7,0)--(7,8);
\draw[dashed, thick] (12,0)--(12,8);
\begin{scope}[shift={(15,0)}]


\draw[fill=black!20!white] (1,0) rectangle (4,8);

\draw[thick,->] (0.7,0)--(4.5,0)node[below right=0.0cm]{};
\draw (2.5,0)node[below =1.cm]  {(b) {time steps}};

\draw[-] (1,0)--(1,-0.3);
\draw[-] (2,0)--(2,-0.3);
\draw[-] (3,0)--(3,-0.3);
\draw[-] (4,0)--(4,-0.3);



\draw (1.5,0) node[below =0.2cm] {$1$};
\draw (2.5,0) node[below =0.2cm] {$2$};
\draw (3.5,0) node[below =0.2cm] {$3$};


\draw[fill=white] (1,0) rectangle (2,2);
\draw (1,0) node[above right =0.2cm] {\small $1$};

\draw[fill=white] (1,2) rectangle (4,4);
\draw (1,2) node[above right =0.2cm] {\small  $2$};

\draw[fill=white] (2,4) rectangle (3,5);
\draw (2,4) node[above right =0.2cm] {\small $3$};

\draw[fill=white] (2,5) rectangle (3,7);
\draw (2,5) node[above right =0.2cm] {\small $4$};

\draw[fill=white] (3,7) rectangle (4,8);
\draw (3,7) node[above right =0.2cm] {\small $5$};

\draw[dashed] (2,0)--(2,8);
\draw[dashed] (3,0)--(3,8);
\end{scope}

\begin{scope}[shift={(21,0)}]

\draw[fill=black!20!white] (1,0) rectangle (4,8);

\draw[thick,->] (0.7,0)--(4.5,0)node[below right=0.0cm]{};
\draw (2.5,0)node[below =1.cm]  {(c) {an optimal solution}};

\draw[-] (1,0)--(1,-0.3);
\draw[-] (2,0)--(2,-0.3);
\draw[-] (3,0)--(3,-0.3);
\draw[-] (4,0)--(4,-0.3);

\draw[-] (0.7,1)--(1,1);
\draw[-] (0.7,2)--(1,2);
\draw[-] (0.7,3)--(1,3);
\draw[ultra thick, -] (0.7,4)--(1,4);
\draw[-] (0.7,5)--(1,5);
\draw[-] (0.7,6)--(1,6);
\draw[-] (0.7,7)--(1,7);
\draw[-] (0.7,8)--(1,8);

\draw (1,2) node[left =0.2cm] {\small  bin 1};
\draw (1,6) node[left =0.2cm] {\small  bin 2};


\draw (1.5,0) node[below =0.2cm] {$1$};
\draw (2.5,0) node[below =0.2cm] {$2$};
\draw (3.5,0) node[below =0.2cm] {$3$};


\draw[fill=white] (1,0) rectangle (2,2);
\draw (1,0) node[above right =0.2cm] {\small $1$};

\draw[fill=white] (1,2) rectangle (4,4);
\draw (1,2) node[above right =0.2cm] {\small  $2$};

\draw[fill=white] (2,4) rectangle (3,5);
\draw (2,4) node[above right =0.2cm] {\small $3$};

\draw[fill=white] (2,0) rectangle (3,2);
\draw (2,0) node[above right =0.2cm] {\small $4$};

\draw[fill=white] (3,0) rectangle (4,1);
\draw (3,0) node[above right =0.2cm] {\small $5$};

\draw[ultra thick] (1,0)--(1,8);
\draw[ultra thick] (4,0)--(4,8);
\draw[ultra thick] (1,4)--(4,4);
\draw[ultra thick] (1,8)--(4,8);
\draw[ultra thick] (1,0)--(4,0);

\draw[dashed] (2,0)--(2,8);
\draw[dashed] (3,0)--(3,8);
\end{scope}

\end{tikzpicture}
\end{center}
\caption{Example of a \TBPP\ instance with 5 items and 3 time steps.}
\label{fig:1}
\end{figure}

In Figure \ref{fig:1}, we give an example to illustrate a \TBPP\ instance with five items having  weights ${w_1=w_2=w_4=2,w_3=w_5=1}$, and bin capacity $W=4$. The starting time and ending time ($s_j$ and $t_j$) are shown in  Figure \ref{fig:1}(a). For instance, item 1 is active in the first two time instants  and item 2 is active {from the second to the last time instant}. {The simultaneously active} sets of  items are : $\overline S_1=\{1\}$, $\overline S_2=\{1,2\}$, $\overline S_3=\{2,3\}$, $\overline S_4=\{2,3,4\}$ and $\overline S_5=\{2,5\}$. The non dominated sets are $\overline S_2$, $\overline S_4$ and $\overline S_5$. In Figure \ref{fig:1}(b) we show the renumbered non-dominated sets
$ S_1=\{1,2\}$, $S_2=\{2,3,4\}$ and $S_3=\{2,5\}$. Finally, Figure \ref{fig:1}(c) reports an optimal {\TBPP} solution using two bins, {the first containing items 1, 2, 4, 5 and the second only item 3.}

One can also use the example to note the difference between \TBPP\ and \2DBPP. Consider an instance of \2DBPP\ defined by the rectangles in Figure \ref{fig:1}, and $3\times 4$ rectangular bins.  An optimal {\2DBPP} solution packs item 3, together with all the other items, in the first bin. {This derives from the fact that in the \2DBPP the rectangles can be shifted not only vertically (as in the \TBPP) but also horizontally (because there are no time windows)}, so a single bin is enough to solve the instance.

We can now show that  the \TBPP\ is a special case of  the \VPP. Given a \TBPP\ instance, define a corresponding \VPP\ instance with the same items and $|T|$ dimensions each with capacity $W$. For each item $j\in N$ define the item weights as
\begin{equation}\label{eq:vpp}
w_{jt}=\left\{\begin{array}{ll}
0 & \text{if\ } t<s_j \text{\ or \ } t\geq t_j\\
w_j& \text{otherwise}\\
\end{array}\right.\quad  {t \in  T}.
\end{equation}
One can see that any \VPP\ solution to this instance is a feasible solution for the \TBPP\ instance, and vice versa, so the \VPP\ generalizes the \TBPP.
\subsection{A polynomial-size model}\label{subsec:poly}
{In this section, we introduce} the first \emph{Integer Linear Programming} (ILP) formulation for the \TBPP.
Let   $I=\{1, 2, \ldots, m\}$  be the set of bins, where $m\leq n$ is an upper bound on the number of bins necessary to {pack} all items.
We introduce a set of binary variables $y$ with the following meaning:
\begin{equation*}
y_{i}=\left\{
\begin{array}{c l}
    1 & \text{if bin $i$ is used,}\\
    0 & \text{otherwise,}
\end{array}\right.
\quad i \in I;
\end{equation*}
and a second set of binary variables $x$ such that:
\begin{equation*}
x_{ij}=\left\{
\begin{array}{c l}
    1 & \text{if item $j$ is packed in bin $i$,}\\
    0 & \text{otherwise,}
\end{array}\right.
\quad i \in I, j \in N.
\end{equation*}
The polynomial-size ILP formulation, called \ILPC\ in the reminder of the paper, reads as follows:
\begin{align}
  && \min\;\; \sum_{i\in I}  y_{i} \label{eq:compact_obj}\\
  && \sum_{i \in I}  x_{ij} & = 1 & j\in N, \label{eq:compact_x=1}\\
  && \sum_{j \in S_t} w_j x_{ij}  &\le W y_{i} &  i\in I, t\in T, \label{eq:compact_capacity}\\
  && x_{ij} &\in \{0,1\}   & i\in I,j\in N, \label{eq:compact_int_x}\\
  && y_{i} &\in \{0,1\}   & i\in I. \label{eq:compact_int_y}
\end{align}
The objective function \eqref{eq:compact_obj} minimizes the number of used bins, constraints \eqref{eq:compact_x=1} impose that each item is packed in one bin, and constraints \eqref{eq:compact_capacity} impose that for each bin and for each time step the total weight of the active items does not exceed the bin capacity.
The optimal solution value of  \ILPC\ is denoted by $z($\ILPC$)$, we use the same notation for all the other models.

\subsection{An exponential-size model}\label{subsec:expo}
{In this section, we} describe a second  formulation for the \TBPP\ characterized by an exponential number of variables associated with all feasible packing patterns, i.e., subsets of items respecting the bin capacity at any time step.
Let $\mathscr{P}$ represent the collection of all feasible packing patterns:
\begin{equation}\label{eq:X}
\mathscr{P} = \left\{P  \subseteq N : \sum_{j \in S_t\cap P} w_j \leq W, \quad t\in T \right\}.\nonumber
\end{equation}
For each  pattern $P\in \mathscr{P}$, we introduce a binary variable $\xi_{P}$ {with the following meaning}:
\begin{equation*}
\xi_{P}=\left\{
\begin{array}{c l}
    1 & \text{if packing pattern $P$ is selected,}\\
    0 & \text{otherwise}.
\end{array}\right.
{\quad P \in \mathscr{P}.}
\end{equation*}
{The exponential-size ILP formulation, called \ILPE\ in the reminder of the paper, reads as follows:}

%
\begin{align}
& &\min\;\;  \sum_{P \in \mathscr{P}}  \xi_{P} \label{eq:extend-obj}\\
&&\sum_{P \in \mathscr{P}: j \in P} \xi_{P} &{=}  1  &   j\in N,  \label{eq:extend-col} \\
&& \xi_{P} &\in \{0,1\}  &  P \in \mathscr{P}. \label{eq:extend-dom}
\end{align}
The objective function \eqref{eq:extend-obj} minimizes the number of packing patterns (bins) used, and constraints \eqref{eq:extend-col} ensure that each item is packed in one bin.
%
%
\section{Lower bounds}\label{sec:lower}
We introduce some lower  bounding techniques which are useful to solve several instances to proven optimality in short computing times.

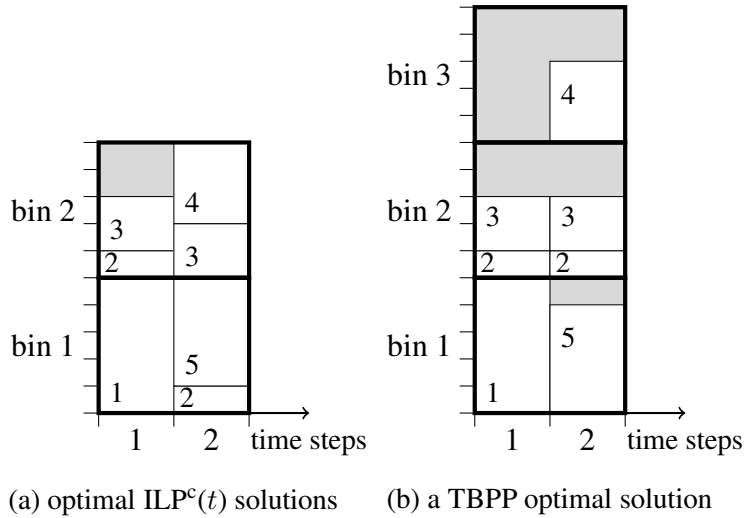
\begin{figure}[ht!]
\begin{center}
\ifx\JPicScale\undefined\def\JPicScale{1}
\unitlength \JPicScale mm
\begin{tikzpicture}[x=\unitlength*10,y=\unitlength*1.8,inner sep=0pt]

\draw[fill=black!15] (0,0) rectangle (2,20);

\draw[fill=white] (0,0) rectangle (1,10); \draw (0,0) node [above right=2mm] {\small 1};
\draw[fill=white] (0,10) rectangle (1,12); \draw (0,10) node [above right=1mm] {\small 2};
\draw[fill=white] (0,12) rectangle (1,16); \draw (0,12) node [above right=2mm] {\small 3};

\draw[fill=white] (1,0) rectangle (2,2); \draw (1,0) node [above right=1mm] {\small 2};
\draw[fill=white] (1,2) rectangle (2,10); \draw (1,2) node [above right=2mm] {\small 5};
\draw[fill=white] (1,10) rectangle (2,14); \draw (1,10) node [above right=2mm] {\small 3};
\draw[fill=white] (1,14) rectangle (2,20); \draw (1,14) node [above right=2mm] {\small 4};

\draw[ultra thick] (0,0) rectangle (2,10);
\draw[ultra thick] (0,10) rectangle (2,20);

\draw (0.5,0) node [below=2mm] {1};
\draw (1.5,0) node [below=2mm] {2};
\draw (0,0)--(0,-1);
\draw (1,0)--(1,-1);
\draw (2,0)--(2,-1);
\foreach \x in {0,2,4,6,8,10,12,14,16,18,20}
{
\draw (0,\x)--(-0.2,\x);
}
\draw[thick,->] (0,0)--(2.8,0) node [below =2mm] {\small time steps};
\draw (1,0) node [below=10mm] {\small (a) optimal ILP\textsuperscript{c}($t$) solutions};

\draw (-0.3,5) node [left] {bin 1};
\draw (-0.3,15) node [left] {bin 2};

\begin{scope} [shift={(5,0)}]
\draw[fill=black!15] (0,0) rectangle (2,30);

\draw[fill=white] (0,0) rectangle (1,10); \draw (0,0) node [above right=2mm] {\small 1};
\draw[fill=white] (0,10) rectangle (1,12); \draw (0,10) node [above right=1mm] {\small 2};
\draw[fill=white] (0,12) rectangle (1,16); \draw (0,13) node [above right=2mm] {\small 3};

\draw[fill=white] (1,0) rectangle (2,8); \draw (1,4) node [above right=2mm] {\small 5};
\draw[fill=white] (1,10) rectangle (2,12); \draw (1,10) node [above right=1mm] {\small 2};
\draw[fill=white] (1,12) rectangle (2,16); \draw (1,13) node [above right=2mm] {\small 3};
\draw[fill=white] (1,20) rectangle (2,26); \draw (1,22) node [above right=2mm] {\small 4};

\draw[ultra thick] (0,0) rectangle (2,10);
\draw[ultra thick] (0,10) rectangle (2,20);
\draw[ultra thick] (0,20) rectangle (2,30);

\draw (0.5,0) node [below=2mm] {1};
\draw (1.5,0) node [below=2mm] {2};
\draw (0,0)--(0,-1);
\draw (1,0)--(1,-1);
\draw (2,0)--(2,-1);

\foreach \x in {0,2,4,6,8,10,12,14,16,18,20,22,24,26,28,30}
{
\draw (0,\x)--(-0.2,\x);
}
\draw[thick,->] (0,0)--(2.8,0) node [below =2mm] {\small time steps};
\draw (1,0) node [below=10mm] {\small (b) a \TBPP\ optimal solution};

\draw (-0.3,5) node [left] {bin 1};
\draw (-0.3,15) node [left] {bin 2};
\draw (-0.3,25) node [left] {bin 3};
\end{scope}

\end{tikzpicture}
\caption{A small example with $LB_0 < z(\ILPC)$}\label{fig:LB0example}
\end{center}
\end{figure}

Observe that for a {single} time step $t \in T$, the \ILPC\ models the  \BPP\ by considering only the items in $ S_t$. We denote this formulation as \ILPC($t$). A valid lower bound for the \TBPP\ is thus:
\begin{equation}
LB_{0}=\max_{t\in T} ~z(\ILPC(t)).
\end{equation}
The following example shows that the optimal solution value of the \TBPP\ may {be strictly greater than} the optimal solution value of the \ILPC($t$) for each $t \in T$.
Consider five items with weights {$w_1=10,w_2=2,w_3=4,w_4=6,w_5=8$}, a bin capacity $W=10$ and two time steps defined by sets $S_1 = \{1,2,3\}$ and $S_2 = \{2,3,4,5\}$. In Figure \ref{fig:LB0example}(a) we show that the optimal  \ILPC($t$) solutions, for $t=1,2$, use two bins, but the optimal \TBPP\ solution  uses  three bins, as depicted in Figure \ref{fig:LB0example}(b). {Two bins can be obtained if and only if items $2$ and item $3$ are packed in different bins in the two time steps. As this does not lead to a feasible TBPP solution, three bins are  necessary for an optimal TBPP solution.}

The difference between $LB_0$ and the optimal \TBPP\ solution value can be greater than one, as shown in the following example. Let $n=24$, $W=3$ and item weights $w_j=2$ ($j=1,\dots,9$),  $w_j=1$ ($j=10,\dots,18$) and $w_j=3$ ($j=19,\dots,24$). Two are the time steps {and the corresponding sets of simultaneously active items are} defined as follows: $S_1 = \{1 \le {j} \le 18\}$ and $S_2 = \{10 \le {j} \le 24\}$.

The optimal \ILPC($t$) solution for $t=1,2$ has $9$ bins as shown in Figure \ref{fig:LB0}(a) (note that, to reduce space,  this figure is drawn with a $90^o$ rotation with respect to the previous ones, i.e., the bins are on the horizontal axis and the time steps in the vertical axis). The optimal \TBPP\ solution has $11$ bins (Figure \ref{fig:LB0}(b)), hence two {units larger} than $LB_0$.
\begin{figure}[ht!]
\begin{center}
\ifx\JPicScale\undefined\def\JPicScale{1}
\unitlength \JPicScale mm
\begin{tikzpicture}[x=\unitlength,y=\unitlength/2,inner sep=0pt]

\begin{scope}[shift={(0,0)}]

\foreach \x in {0,1,2,3}
{
\draw (-3,\x*10)--(0,\x*10);
}
\draw[thick,->] (0,0)--(0,35);

\foreach \x in {1,2,3,4,5,6,7,8,9}
{
\draw[fill=white] ({(\x-1)*10},0) rectangle ({(\x-1)*10+10},20);
\draw[fill=white] ({(\x-1)*10},20) rectangle ({(\x-1)*10+10},30);
\draw ({(\x-1)*10+5},-5) node {\footnotesize bin \x};
\draw ( {(\x-1)*10+5},5) node {\footnotesize \x};
\draw[ultra thick] ({(\x-1)*10},0) rectangle ({(\x-1)*10+10},30);
}
\foreach \x in {10,11,12,13,14,15,16,17,18}
{
\draw ( {(\x-10)*10+5},25) node {\footnotesize \x};
}

\end{scope}

\begin{scope}[shift={(0,-50)}]

\foreach \x in {0,1,2,3}
{
\draw (-3,\x*10)--(0,\x*10);
}
\draw[thick,->] (0,0)--(0,35);

\foreach \x in {1,2,3,4,5,6}
{
\draw[fill=white] ({(\x-1)*10},0) rectangle ({(\x-1)*10+10},30);
\draw ({(\x-1)*10+5},-5) node {\footnotesize bin \x};
\draw[ultra thick] ({(\x-1)*10},0) rectangle ({(\x-1)*10+10},30);
}
\foreach \x in {19,20,21,22,23,24}
{
\draw ( {(\x-19)*10+5},15) node {\footnotesize \x};
}
\foreach \x in {7,8,9}
{
\draw[fill=white] (\x*10-10,0) rectangle (\x*10,10);
\draw[fill=white] (\x*10-10,10) rectangle (\x*10,20);
\draw[fill=white] (\x*10-10,20) rectangle (\x*10,30);
\draw (\x*10-5,-5) node {\footnotesize bin \x};
\draw[ultra thick] ({\x*10-10},0) rectangle ({\x*10},30);
}
\foreach \x in {16,17,18}
{
\draw ( {(\x-16)*10+65},25) node {\footnotesize \x};
}
\foreach \x in {13,14,15}
{
\draw ( {(\x-13)*10+65},15) node {\footnotesize \x};
}
\foreach \x in {10,11,12}
{
\draw ( {(\x-10)*10+65},5) node {\footnotesize \x};
}

\draw (0,-15) node [below right=5mm] {\footnotesize (a) Optimal \ILPC($t$) solutions for time steps 1 (above) and 2 (below)};

\end{scope}


\begin{scope} [shift={(0,-125)}]

\foreach \x in {0,1,2,3}
{
\draw (-3,\x*10)--(0,\x*10);
}
\draw[thick,->] (0,0)--(0,35);

\foreach \x in {1,2,3,4,5,6,7,8,9}
{
\draw[fill=white] ({(\x-1)*10},0) rectangle ({(\x-1)*10+10},20);
\draw ({(\x-1)*10+5},-5) node {\footnotesize bin \x};
\draw ( {(\x-1)*10+5},5) node {\footnotesize \x};
}

\foreach \x in {16,17,18}
{
\draw[fill=white] ({(\x-10)*10},20) rectangle ({(\x-10)*10+10},30);
\draw ( {(\x-10)*10+5},25) node {\footnotesize \x};
}

\foreach \x in {10,11}
{
\draw[fill=white] (\x*10-10,0) rectangle (\x*10,10);
\draw[fill=white] (\x*10-10,10) rectangle (\x*10,20);
\draw[fill=white] (\x*10-10,20) rectangle (\x*10,30);
\draw (\x*10-5,-5) node {\footnotesize bin \x};
}
\foreach \x in {10,11,12}
{
\draw ( 95,{(\x-10)*10+5}) node {\footnotesize \x};
}
\foreach \x in {13,14,15}
{
\draw ( 105,{(\x-13)*10+5}) node {\footnotesize \x};
}
\foreach \x in {1,2,3,4,5,6}
{
\draw[fill=black!15] ({(\x-1)*10},20) rectangle ({(\x-1)*10+10},30);
}

\foreach \x in {1,2,3,4,5,6,7,8,9,10,11}
{
\draw[ultra thick] ({(\x-1)*10},0) rectangle ({(\x-1)*10+10},30);
}

\end{scope}

\begin{scope}[shift={(0,-175)}]

\foreach \x in {0,1,2,3}
{
\draw (-3,\x*10)--(0,\x*10);
}
\draw[thick,->] (0,0)--(0,35);

\foreach \x in {1,2,3,4,5,6}
{
\draw[fill=white] ({(\x-1)*10},0) rectangle ({(\x-1)*10+10},30);
\draw ({(\x-1)*10+5},-5) node {\footnotesize bin \x};
}
\foreach \x in {19,20,21,22,23,24}
{
\draw ( {(\x-19)*10+5},15) node {\footnotesize \x};
}
\foreach \x in {7,8,9}
{
\draw[fill=white] (\x*10-10,0) rectangle (\x*10,10);
\draw (\x*10-5,-5) node {\footnotesize bin \x};
}
\foreach \x in {16,17,18}
{
\draw ( {(\x-16)*10+65},5) node {\footnotesize \x};
}
\foreach \x in {10,11}
{
\draw[fill=white] (\x*10-10,0) rectangle (\x*10,10);
\draw[fill=white] (\x*10-10,10) rectangle (\x*10,20);
\draw[fill=white] (\x*10-10,20) rectangle (\x*10,30);
\draw (\x*10-5,-5) node {\footnotesize bin \x};
}
\foreach \x in {10,11,12}
{
\draw ( 95,{(\x-10)*10+5}) node {\footnotesize \x};
}
\foreach \x in {13,14,15}
{
\draw ( 105,{(\x-13)*10+5}) node {\footnotesize \x};
}
\foreach \x in {7,8,9}
{
\draw[fill=black!15] ({(\x-1)*10},10) rectangle ({(\x-1)*10+10},30);
}

\foreach \x in {1,2,3,4,5,6,7,8,9,10,11}
{
\draw[ultra thick] ({(\x-1)*10},0) rectangle ({(\x-1)*10+10},30);
}

\draw (0,-15) node [below right=5mm] {\footnotesize (b) Optimal \TBPP\ solution: time steps 1 (above) and 2 (below) };

\end{scope}

\end{tikzpicture}
\caption{Example showing that $LB_0$ may have a gap larger than one.}\label{fig:LB0}
\end{center}
\end{figure}

Let us now consider the continuous relaxation of  \ILPC\ and \ILPC($t$),  denoted, respectively, by  \LPC\ and \LPC($t$).
A second valid \TBPP\ lower bound is  given by
$z (\LP\textsuperscript{c}).$
For {any} $t \in T$, it is known (see e.g. \citet{MT90}) that the optimal solution value of \LPC($t$)   can be computed as:
%
$z ({\rm LP}\textsuperscript{c}(t)) = \frac{1}{W}\sum_{j\in S_t}w_j,$
%
from which we immediately have
\[
z(\LPC)=\max\{z(\LPC(t)), t\in {T}\}=\frac{1}{W}\max\left\{\sum_{j\in S_t}w_j, t\in {T}\right\}.
\]
Given an instance of the \TBPP\, let $\tilde{t} =\argmax\{\sum_{j\in S_t}w_j: t\in T \}$, then a lower bound from the continuous relaxation of \ILPC\ is
\begin{equation}\label{LB1_T}
LB_{1}= \lceil z (\LPC(\tilde{t})) \rceil.
\end{equation}
 A straightforward implementation of $LB_1$ requires an $O(n^2)$ computing time, but this bound can be improved as shown in the next property.\\

\begin{property}\label{pro:LB1}
Given an instance of the \TBPP, lower bound $LB_1$ can be computed in $O(n\log n)$ time.
\end{property}
\begin{proof}
Remind that the items are sorted by non-decreasing starting times. Compute a list $\cal L$ containing the items sorted by non-decreasing ending time $t_j$.
Start with a total item weight ${\cal W}=0$ and consider the items one at a time, in the original order. For each item $j$, add $w_j$ to $\cal W$ and remove from list $\cal L$ each item $k$  with $t_k\leq s_j$, by reducing $\cal W$ of the corresponding weights $w_k$. {Store the maximum value of  $\cal W$ obtained} and let $\tilde{t}$ denote the corresponding time step. The initial sorting can be implemented in  $O(n\log n)$, while the scanning of the items requires $O(n)$ because each item is considered exactly once in each list, and the thesis follows.
\end{proof}

 Another lower bound, which requires an intermediate computing effort between that of $LB_0$ and $LB_1$, can be obtained by solving a single \BPP\ on  time step $\tilde{t}$, thus obtaining: $$LB_0(\tilde{t}) = ~z(\ILPC(\tilde{t})).$$

\subsection{Preprocessing with item weight lifting}\label{subsec:lifting}

In order to improve lower bound {$LB_1$, see \eqref{LB1_T},} we can apply lifting techniques  (see, e.g., \citet{DDI12}) which try to increase the item weights, while ensuring that the lifted instance has the same optimal solution {value} of the original one.

Given an item $j\in N$, let $\gamma(j) = \{k\in N\setminus\{j\}: \exists t\in T \text{\ such that\ } j,k\in S_t \}$ denote the set of items that are active in at least one time step where $j$ is active. We can use the following ILP model to compute the maximum possible loading of {a bin} where {item} $j$ is packed:
\begin{equation}
 \sigma(j) = \max  \left\{\sum_{k \in \gamma(j)} w_k x_k  \label{eq:SS_1} :
\sum_{k \in \gamma(j)} w_k x_k  \le W - w_{j}, x_k  \in \{0,1\}, k \in \gamma(j)\right\}.
\end{equation}
Given the optimal value $\sigma(j)$ {(or any valid upper bound on $\sigma(j)$)} we know that the bin where $j$ is packed, in the time steps where $j$ is active, has at least  an empty space $W- w_{j} - \sigma(j)$, so we can lift {the weight of $j$ as:}
 \begin{equation}\label{eq:lift1}
w_j = w_{j} + (W- w_{j} - \sigma(j)) = W - \sigma(j) \tag{Lift-1}.
 \end{equation}
Model \eqref{eq:SS_1} represents a Subset Sum Problem (\SSP) (see, e.g., \citet{MT90}).
To implement this lifting procedure we have to execute a  first step which solves $n$ {\SSP}s to try to lift the {item weights}.
{Once an item $j$ has been lifted, a new iteration can be executed to try to further lift the items in $\gamma(j)$.}

Another lifting procedure can be obtained by considering each  time step $t\in T$, such that $j\in S_t$, and defining the following \SSP:
\begin{equation}\label{eq:SS_2}
 \sigma(j,t) = \max  \left\{\sum_{k \in {S}_t\setminus\{j\}} w_k x_k  : \sum_{k \in {S}_t\setminus\{j\}} w_k x_k  \le W - w_{j}, x_k \in \{0,1\},  k \in {S}_t\setminus\{j\} \right\}.
\end{equation}
Similarly to the previous case, we see that, for the given time step $t$, we can lift the weight of the item $j$ to $  W - \sigma(j,t)$. By considering all the time steps in which $j$ is active, we define the valid lifting:
\begin{equation}\label{eq:lift2}
w_{j} = n_{t \in T: j\in S_t } (W-\sigma(j,t)).\tag{Lift-2}
\end{equation}
The implementation of this lifting requires to solve $n\times |T|$ {\SSP}s for each iteration.\\
\begin{property}\label{obs:lift2}
For a given item $j\in N$, lifting \text{Lift-2} dominates lifting \text{Lift-1}.
\end{property}
\begin{proof}
Let $\tilde{t}$ denote the time instant giving the minimum value in \eqref{eq:lift2}. The thesis immediately follows because $S_{\tilde{t}}\subseteq \gamma(j)$.
\end{proof}
We note that, although observation \ref{obs:lift2} guarantees that Lift-2 dominates Lift-1 for a given item, it does not guarantee that the total lift on all the items {obtained by Lift-2 dominates that obtained} by Lift-1.

In the following, we will denote as $LB^I_2$ and $LB^{II}_2$ the value of $LB_1$ computed with weights lifted using the first and second method, respectively, and a single lifting iteration.
\subsection{Lower bound from the exponential-size formulation}\label{subsec:exp-lb}
Another lower bound can be obtained by computing the continuous relaxation of the exponential-size formulation 
\ILPE, that we denote as \LPE\ in what follows. We have thus
\begin{equation}\label{eq:LB3}
LB_3 = {\lceil z(\LPE) \rceil}.
\end{equation}
{ The following property states that the quality of the lower bound obtained solving the LP relaxation of \ILPC\ is dominated by its counterpart associated with \ILPE.\\}

\begin{property}\label{obs:obs2}
The optimal value of \LPE\ is greater than or equal to the optimal value of \LPC.
\begin{proof}
We start the proof by showing that any feasible solution of \LPE\ can be transformed into a feasible solution of \LPC\ preserving  its objective function value.
W.l.o.g., we consider the case in which $m=n$. {Any optimal basic solution of \LPE\  can have at most $n$ non-zero variables}; let $\xi^*$ denote a feasible solution of \LPE, and
$P(i)$ be a function that returns the {$i$-th} active pattern in $\xi^*$ (i.e., $P(i)$ is the {$i$-th} pattern associated to a strictly positive variable).
We can define a solution $(x^*,y^*)$ of \LPC\ as follows: for each bin $i \in I$ and for each $j \in N$, we can set:
\begin{equation*}
  y^*_i =  \xi^*_{P(i)}
  \quad \text{and } \quad
  x^*_{ij} = \xi^*_{P(i)}.
\end{equation*}
 In case $ \sum_{i \in I} x_{ij}^*>1$, for some $j$, we can arbitrarily reduce the $x_{ij}^*$ in order to sum up 1. By construction,  the solution $(x^*,y^*)$ is feasible for \LPC\ and has the same objective function value.

We then show a case where the optimal value of \LPE\ is strictly larger than the optimal value of \LPC.
Consider the instance presented in Figure \ref{fig:LB0example} with 5 items and 2 time steps. The optimal solution value of   \LPC\ is equal to 2. {This optimal  value can be obtained for instance by the solution $y_1^*=y_2^*=1$ and $x_{ij}=\frac{1}{2}$ for $i=1,\dots,5$ and $j=1,2$ (all the other \LPC\ variables are set to zero).}
An optimal solution for \LPE\ is defined by the four feasible packing patters $P_1=\{2,3\}$, $P_2=\{2,4\}$, $P_3=\{1,5\}$, $P_4=\{3,4\}$, and by the corresponding variables $\xi_{P_1}=\xi_{P_2}=\xi_{P_4}=0.5$ and $\xi_{P_3}=1$. The optimal solution value is  $z(\LPE)=2.5>z(\LPC){=2}$.
%
%
\end{proof}
\end{property}
\section{Upper bounds}\label{sec:ub}
We start by  describing some simple heuristic algorithms that can be used to provide  approximate solutions in short computing times.

\paragraph{{Greedy algorithm}}
{The first method we introduce is based on a sequence of greedy algorithms} derived from the well known \emph{First-Fit} algorithm for \BPP\ (see, e.g., \citet{J85}).
The First-Fit performs $n$ iterations by considering one  item at a time, in a given order. At the beginning no bin is open (used). In the first iteration, a single bin is opened and the first item is packed in it. At each iteration $j> 1$, the algorithm looks for the first open bin in which $j$ fits,  if any. If this bin exists $j$ is packed in it, otherwise a new empty bin is opened and used to pack $j$. The algorithm runs in $O(n^2)$ time.

We extend the First-Fit for the \BPP\ to the \TBPP, by simply checking, for each item $j$ and tentative bin $i$, if packing $j$ in $i$ is feasible for all the time steps in $[s_j,t_j)$. The time complexity increases to $O(n^2|T|)$, but in practice these methods are very fast. We implemented two versions: in First-Fit-1 we do not make any sorting of the items, but we consider them in their natural order given by increasing $s_j$ values. In First-Fit-2 we sort the items by non-decreasing $w_j$ values. The two versions were run on the original instances and on those lifted with Lift-1 (see Section \ref{subsec:lifting}). In the following, we call ${UB_1}$ the best solution value obtained by running First-Fit-1 and First-Fit-2 on the original and lifted instances.

\paragraph{Rolling horizon heuristic}
We developed a more powerful heuristic  using a rolling horizon concept. Let $\Delta\in [1,\dots,|T|]$ denote the width of a heuristic time window.
{Our \emph{rolling horizon heuristic}} performs $\lceil\ {|T|}/{\Delta}\ \rceil$ iterations, and at each iteration it solves  model \ILPC\ with a limited number of time-step constraints. More specifically, let $\ILPC_{k\Delta}$, with $k=1,\dots, \lceil\ {|T|}/{\Delta}\ \rceil$, denote the restricted \ILPC\ model obtained by:
(i) defining constraints \eqref{eq:compact_capacity} only for time steps in $T_k=[(k-1)\Delta+1, n(k\Delta, |T|)]$, and
(ii) selecting only the $x$ variables associated with the items in $N_k=\cup_{t\in T_k} S_j$ (i.e., the  active items in $T_k$).
 If a feasible solution for $\ILPC_{k\Delta}$ is found, we pack the items in $N_k$ according to this solution, by fixing, the corresponding $x$ variables to one.
 These items are no longer considered for possible reassignment. Due to the variable fixing the next model $\ILPC_{(k+1)\Delta}$ starts from the packing of the items in $\cup_{h=1}^k N_h$ and packs the new items from $N_{k+1}$ in the residual space of the opened bins, or in new bins.
 If for some iteration $k$ no feasible solution to $\ILPC_{k\Delta}$ is found, the algorithm terminates with no \TBPP\ solution.
 We will call $UB_2$ the solution value obtained by the rolling horizon heuristic.

A further upper bound $UB_3$ has been developed using a truncated version of the exact branch-and-price algorithm of Section \ref{sec:BP}. We will give details on this heuristic in that section.
\section{Solving the Exponential-Size formulation}\label{sec:BP}
In Section \ref{sec:Models}, we introduced a polynomial-size formulation \ILPC\ and an exponential-size formulation \ILPE. Formulation \ILPC\ can be
explicitly written, also for large size instances, and solved using a generic ILP solver. Formulation
\ILPE, instead,  has exponentially many variables that cannot be explicitly enumerated for large-size instances.
\emph{Column Generation} (CG) techniques are then necessary to efficiently solve the continuous relaxation of \ILPE\ (we  refer the interested reader to, e.g., \citet{desaulniers2006column} for further details on column generation).
In the following, we present a new branch-and-price framework for \ILPE.
Two are the main ingredients of a branch-and-price algorithm: (i) a column generation algorithm to solve the Linear Programming Relaxation of the exponential-size integer model, and (ii) a branching scheme. We discuss separately these two aspects in the next sections, before  introducing a heuristic based on this framework, and our final overall algorithm for the \TBPP.
\subsection{Solving the Linear Programming Relaxation of \ILPE}\label{sub:solvingLPe}
Model \LPE,  initialized with a subset of variables (columns) defining feasible solutions, is called the \emph{Restricted Master Problem} (RMP).
{In our implementation, we initialize the RMP with the columns associated with the best solution provided by the Greedy algorithm of Section \ref{sec:ub}.}
Its solution provides a (sub)-optimal primal solution. To find the optimal solution we consider the
 dual of this LP:
%
{
\begin{align}
\max  \left\{\sum_{j \in N}  \pi_{j} : \quad \sum_{j \in P} \pi_j \le 1,\; P \in \mathscr{P}, \quad \pi_{j} \geq 0,\;  j \in N \right\},
\end{align}}
where $\pi_j$  is the dual variable associated with the $j$-th constraint \eqref{eq:extend-col}.
{
It is worth noticing that we use the `$\geq$' sign in constraint \eqref{eq:extend-col}, instead of the `=' sign, to have non negative dual variables. This choice does not change the problem, since any solution of \eqref{eq:extend-obj}-\eqref{eq:extend-dom} which packs an item in more than one bin can be transformed into a \TBPP\ solution by arbitrarily removing the item from all the used  patterns, but one.
}
A violated dual constraint induces a negative reduced cost in the primal problem, so the corresponding primal variable must be added
to the RMP to find an optimal solution.
 Accordingly, the column generation performs a number of iterations where violated dual constraints are added to the RMP in form of primal variables, and the RMP is re-optimized, until no violated dual constraint exists. At each iteration, the so-called \emph{Pricing Problem} is solved. This problem asks to determine (if any) a packing pattern $P^* \in \mathscr{P}$ for which the associated
dual constraint is violated, i.e., such that
\begin{equation}
    \sum_{j \in P^*} \pi^*_j > 1, \label{pricing_cond}
\end{equation}
where $\pi^*$ is the optimal vector of dual variables for the current RMP.

If a packing pattern $P^*$ has dual weight larger than one (that is, the reduced cost is negative), the associated column is added to the RMP and the problem is re-optimized. If, on the other hand, the dual weight is not larger than 1, by linear programming optimality conditions no column can improve the objective function of the RMP and therefore \LPE\ is solved to optimality.

In the solution of the pricing problem  we want to find a violating pattern $P^*$, or to prove that no one exists. The following
 0-1 \emph{Temporal Knapsack Problem} (TKP) models the separation using
 a binary variable $z_j$ that takes value 1 if and only if item $j\in N$ is selected in subset $P^*$:
{
\begin{align}
 z_s(\pi^*)=\max \left\{\sum_{j \in N}\pi^*_j z_j
:\quad \sum_{j \in S_t} w_j z_j \le W,\; t\in T, \quad z_{j} \in \{0,1\}, \; j\in N \right\}. \label{eq:slave_3}
\end{align}
}
If $z_s(\pi^*) >1 $ a violating pattern $P^*=\{j\in N: z_j=1\}$ has been found, otherwise the RMP solution is optimal.

\subsection{Branching schemes for \ILPE}
The design of a branching scheme is crucial for the performance of a branch-and-price algorithm (see, e.g.,  \citet{Vanderbeck2011}).
In the following, we describe two branching scheme adopted in our new branch-and-price framework. Two are the main properties that a branching rule should hold.
Firstly, it is a complete scheme, i.e., it ensures that integrality can be imposed in all cases.
Secondly, it  does not require modifications to the master problem and it does not impact much on the pricing algorithm.
The latter property means that an ideal branching does not alter the structure of the  pricing problem so that the same algorithm can be applied during the entire search.
In the following, we denote with  $\xi^*$ a fractional solution  to \LPE at a given node of the branching tree, and with $\mathscr{\widehat P} \subseteq \mathscr{P}$ the set of columns in the RMP at the node.

\paragraph{Branching BR-1} This is the standard  branching rule, which  selects a  variable $\xi^*_P$ with fractional value and separates the current node into two children nodes, by imposing, respectively, $\xi^*_P=1$ and $\xi^*_P=0$. The advantage of this rule is to force the algorithm to  find  feasible solutions in short time, in particular when a deep-first exploration rule is used and the left branch  (with $\xi^*_P=1$) is selected first.
The RMP can easily incorporate the branching constrains, because one can directly impose the branching variable value. The pricing on the left branch is easy, because we can just remove the items in $P$ from the pricing problem, but on the right branch it is necessary to add a cut to avoid generating $P$ again. {To achieve this,} we use the cut
\begin{equation}
\sum_{j\in N\setminus P}z_j -\sum_{j\in  P}z_j \geq 1-|P|.\label{eq:branch1-cut}
\end{equation}
\paragraph{Branching BR-2}
The second branching strategy we propose is inspired by the \textit{Rayan-Foster branching} scheme and it  preserves the pricing algorithm in part of the branching nodes.
This  rule is designed to impose that each couple of items are either packed in the same bin or in different ones.
This rule has been proposed for {branch-and-price algorithms based on set-covering formulations (see, e.g.,  \citet{Barnhart98})} and used to derive several effective exact algorithms for the \BPP, see e.g., {\citet{Wei2018}}.
A couple of items  $r$ and $s$   $\in N$  is {\it fractionally packed} if:
\begin{align}
 \sum_{P \in \mathscr{\widehat P} :  r,s \in P} \xi_{P}^* = \gamma,\label{eq:fractional}
\end{align}
with $\gamma$ fractional.
 In case more than one pair of fractionally packed items exist, we select the couple  $r$ and $s$ associated with the first fractional $\gamma$ value.  Two children nodes are then created:
\begin{itemize}
    \item in the first node we force items $r$ and $s$ to be packed in the same bin;
    \item in the second node we force items $r$ and $s$ to be packed in different bins.
\end{itemize}
 This  branching rule can be implemented without any additional constraint to  the RMP, indeed it is enough to remove from RMP the columns that does not respect the  rule.
For the pricing problem we can see that in the left branch we can force  a couple $(r,s)$ of items to be packed together  by   replacing the couple with a \textit{super item}, say $\hat \jmath$ with $s_{\hat \jmath}=n(s_r,s_s)$, $t_{\hat \jmath}=\max(t_r,t_s)$ and weight depending on the time steps:
\[
w_{\hat\jmath t}=\left\{\begin{array}{ll}
  w_r & \text{if\ } r\in S_t,\ s\not\in S_t,\\
  w_s & \text{if\ } r\not\in S_t,\ s\in S_t,\\
  w_r+w_s & \text{if\ } r\in S_t,\ s\in S_t,
\end{array}\right.
\quad t\in   [s_{\hat \jmath}, t_{\hat \jmath}].
\]
{The pricing problem is modified by substituting {the capacity constraints in  \eqref {eq:slave_3}} with  constraints $\sum_{j \in S_t} w_{jt} z_j \le W, t\in T$, so the problem remains a \TKP\ but  {with time step dependent item weights}.}

In the right branch we  enforce the two items to belong to different patterns. In this case  the pricing problem {\eqref {eq:slave_3}}
must be changed by introducing the so called \emph{conflicts} between  items. Unfortunately the conflicts cannot be directly managed with some modification of the input instance,  and we are forced to solve the problem as a generic MIP with the additional constraint $z_r+z_s\leq 1$.

The branching rule BR-1 is clearly complete because we have no variable to branch only when all values are integer, and the solution is integer. The following observation states that also the  branching rule BR-2 is complete for \ILPE:\\
\begin{property}\label{obs:obs3}
    The  branching rule BR-2 provides a complete branching scheme for model \ILPE.
\begin{proof}
If at each node of the branching tree we select a pair of items and  fix them to be packed in the same bin or in separate bins, then, after at most $O(n^2)$ branchings, all  couples have been fixed and all solutions enumerated. To prove the thesis, it is enough to show that we can always find two items providing a fractional $\gamma$ value,  as defined in \eqref{eq:fractional}.

In \citet{Barnhart98}, it is proved that for any 0-1 constraint matrix $A$ (as for the case of \LPE), if a basic solution $\xi^*$ to $A \xi = 1$ is fractional, then there exist two rows $r$ and $s$ such that:
        \begin{equation}
            0 < \sum_{P \in \mathscr{P} \ : \ r, s \in P} \xi_{P}^* < 1. \label{eq:barnhart}
        \end{equation}
As we adopted a covering formulation ($A\xi \geq 1$), we can use the above  property to prove that the required items exist if $\xi^*$ determines  at least two strict constraints (i.e., it is satisfied with the ``=" sign). We show that these two rows always exists by analyzing all possible cases.
We first note that in an optimal solution there must be at least one strict constraint, otherwise it is possible to improve the solution by reducing the value of one variable, until the l.h.s. of one row takes the value one. If there is at least another strict constraint we are done, otherwise let $r$ denote the row of the unique strict constraint.
 If there is  a packing $P$ not containing item $r$ and such that $\xi^*_P > 0$, we can improve the solution by reducing $\xi^*_P$ until a second constraint is strict, and this case is closed. Otherwise  let  $\mathscr{\widehat P^+}$ denote the subset of columns of  $\mathscr{\widehat P}$ with a positive $\xi^*$ value. In the remaining case, $r\in  P$ for all $P\in \mathscr{\widehat P^+}$, and for any $s\in N\setminus\{r\}$ equation \eqref{eq:fractional} defines an integer $\gamma$ value (otherwise we have found the $r,s$ couple to be used for branching). But the latter case is not possible because it implies that any item is packed within $r$ in all  columns of $\mathscr{\widehat P^+}$, which results to be identical.
     \end{proof}
\end{property}
\subsection{Exponential-size formulation based heuristic}\label{subsec:exp-heur}
Effective heuristic algorithms can be obtained by contaminating a branch-and-price framework with {one or more heuristic rules} that cut parts of the search tree. The \emph{diving} metaphor in an LP-based branch-and-bound tree foresees a search that plunges deep into the enumeration
tree by selecting a branch with some heuristic rule at each node.

In our case, we have implemented the branch-and-price with branching rule BR-1 (which fixes a variable $\xi_P$ at a time to 1 and 0, respectively), but we have limited the number of possible branching using a \emph{token-based} rule. Fixing a variable at value 1 (i.e., fixing the packing of a bin) reduces the problem size and drives the algorithm to find a feasible solution in a short time. For this reason, we allow the algorithm to perform all left branches ( which fix the variable to 1), but we restrict the number of right  branches. More specifically, at the beginning the algorithm has $K$ token available for right branching. When a tree node $k$ is separated using a right branch, the number of available tokens is reduced by one. When no more  token  exist, the algorithm can perform only left branches. If backtracking occurs and the search return to node $k$ the token possibly used for a right branch returns available and the total number of tokens is increased by one. One can see that with $K\geq 1$ tokens at most $O(n2^{(K-1)})$ solutions are generated, while a single solution is generated for $K=0$.
\subsection{Overall algorithm}\label{subsec:exp-overall}
The final algorithm we propose to solve the \TBPP\ is based on a combination of the above approaches. In particular, in a first phase,
we start by computing the lower bound $LB_0$ and an upper bound, say $UB_0$, obtained by running the Greedy and the H-Rolling heuristics, and choosing the best solution. If $UB_0=LB_0$, then  we stop. Otherwise, in the second phase,  we compute the continuous relaxation \LPE. We check again if $UB_0$ provides the optimal solution by comparing it with $LB_3$, and stopping if equality holds. In the third phase, we run the tree-exploration of the H-Diving heuristic, starting from the  continuous relaxation \LPE\ already computed while  evaluating $LB_3$.
Let  $UB_3$ be the corresponding heuristic solution value.  If $UB_3=LB_3$ we are done, otherwise we proceed to the final phase running the branch-and-price method of Section \ref{sec:BP}, using branching rule BR-2. The algorithm, called B\&P$^+$ in the following, is resumed in Algorithm \ref{alg-B&P}.
We observe that phase 2 is one of the most time consuming, because it requires to solve \LPE\ using the column generation approach. Once \LPE\ is solved the next tree-search performed by H-Diving with a single token, is usually fast. In phase 4, instead, the tree-search  using branching rule BR-2 without restrictions, may use large computing times. See Section \ref{sec:computational} for details.

\begin{algorithm}[bht]
\caption{{{B\&P$^+$}}}\label{alg-B&P}
\begin{algorithmic}[1]
\State {Compute $L=\max(LB_0,LB_1,LB_2)$}; \Comment Phase 1
\State {Run First-Fit-1 and First-Fit-2 on the original and lifted instance, run H-Rolling, and let $U$ be the minimum solution value};
\State {\textbf{if} $U=L$ the solution is optimal \textbf{then} \textbf{return};}
\State {Compute $L=\max(L,LB_3)$}; \textbf{if} $U=L$ the solution is optimal \textbf{then} \textbf{return}; \Comment Phase 2
\State {Run H-Diving with one token and possibly improve $U$;} \Comment Phase 3
\State \textbf{if} $U=L$ the solution is optimal \textbf{then} \textbf{return};
\State {Run the branch-and-price with branching rule BR-2, possibly improving $L$ and $U$ \Comment Phase 3}
\State {\textbf{return}}.
\end{algorithmic}
\end{algorithm}

\section{Computational Results}\label{sec:computational}

The experiments have been performed on a computer with a 3.10 GHz 4-core Intel Xeon processor and 16Gb RAM, running a 64 bits Ubuntu Linux operating system version {14.04.5}. The algorithms were coded in \texttt{C++} and all the codes were compiled with \texttt{gcc 6.2} and \texttt{-O3} optimizations.
To solve the linear relaxations of the models and the required {ILPs}, we used  \texttt{Cplex 12.7}, run on a single-thread (parameter \texttt{CPX\_PARAM\_THREADS} set to one).

We build our testbed starting from the \TKP\ instances called ``I" in \citet{CFM13}.
The testbed ``I" consists of one hundred instances which are further divided into ten classes generated using different values of some inout parameters (see Table \ref{tab:param}), and the following rules, which allows to generate only the {inclusion-wise non-dominated} sets {of simultaneously active items} (see \citet{CFM13} for further details).
Given the size $|T^o|$, for each $t\in \{1,\dots,|T^o|\}$ the number of
tasks in $S_t$  is uniformly distributed in $[a_{n},a_{\max}]$. If $t>1$ then $\beta$ percent tasks from $S_{t-1}$ are randomly selected and inserted in $S_t$, where $\beta$ is  uniformly distributed in $[b_{n},b_{\max}]$.
Item weights  are uniform random integer from $[10,100]$. The bin capacity is $W=100$ for all instances.
 For Classes I-IV, $|T^o|=2688+128(i-1)$  ($i=1,\ldots,10$). For Classes V-X, $|T^o|=768+128(i-1)$
($i=1,\ldots,10$).

\begin{table}[ht!] \caption{Parameter values used to generate the test  instances\label{tab:param}} \centering
\footnotesize \tabcolsep 3.6pt
\renewcommand
\arraystretch{1.2}
\begin{tabular}{cccccccccccccc}
\hline
Class &   $a_{n}$  &   $a_{\max}$  &   $b_{n}$  &   $b_{\max}$ &&Class &   $a_{n}$  &   $a_{\max}$  &   $b_{n}$  &   $b_{\max}$   \\
\hline \hline
I   &   10  &   10  &   90  &   95   && VI  &   30  &   30  &   70  &   90   &      \\
II  &   15  &   15  &   90  &   95   && VII &   30  &   30  &   90  &   95   &      \\
III &   20  &   20  &   90  &   95   && VIII&   25  &   35  &   90  &   95   &      \\
IV  &   25  &   25  &   90  &   95   && IX  &   25  &   35  &   70  &   90   &      \\
V   &   30  &   30  &   90  &   95   && X   &   30  &   40  &   90  &   95   &      \\
\hline
\end{tabular}
\end{table}

Due to the fact that the \TBPP\ instances share the same data structure with those of the \TKP,
except for the absence of item profits, these instances could be directly used to define \TBPP\ instances. However the \BPP\ is normally much harder to solve than a KP, and the same happens with the time versions here considered. Therefore we reduced the size of the original instances as follows. Let $T^o$ denote the time steps in an original instance, we build \TBPP\ instances with $|T|\in \{5,10,15,20,30,\dots,150\}$ by  extracting from the original instance all the items in $\cup_{t=1}^{|T|} S_t$, and disregarding the profits.

Instances of classes  VI and IX have $b_{n}$ {smaller} than the other classes.
This implies that when generating a set $S_t$ a smaller number of tasks are taken from $S_{t-1}$, and thus a larger number of new tasks exist in $S_t$.
This fact, however, does not translate into a greater difficulty of the class, as will be seen in the remaining of the section.

\subsection{Lower bounds}\label{subsec:comp-lb}
Our first lower bound $LB_0$ (see Section \ref{sec:lower}) requires to solve $|T|$ \BPP\ instances. To solve the {\BPP}s we implemented  the arc-flow algorithm proposed in {\citet{V02}}. Lower bound $LB_1$ can be computed in $O(n\log n)$ as shown by Property \ref{pro:LB1}, while $LB_2^I$ requires to solve $n$ Subset Sum problems for each iteration (see Section \ref{subsec:lifting}). We executed a single iteration by solving the \SSP\ with \texttt{Cplex}. The next lower bound $LB_2^{II}$ is still based on a lifting procedure, and requires the solution of $O(n|T|)$ {\SSP}s for each iteration. Again we used a single iteration and \texttt{Cplex} for solving the subproblems.

Table \ref{tab:LBsALL-timeSteps} reports the results of the lower bounds for instances with $|T|=(10,20,\dots,100)$. Each row corresponds to 100 instances from the ten classes (ten instances for each class). The column labeled  $|T|$ gives the number of time steps in the 100 instances, the column labeled `$|N|$' gives the average number of items. Columns `\# opt', `avg gap' and `max gap',  report, respectively,  the number of times the bound is equal to the optimal solution value and  the average and  maximum absolute gap with respect to the optimum. (The optimal solution values are computed with the exact methods evaluated in the next Section \ref{subsec:comp-exact}, possibly running them for long times, when no method is able to produce a proven optimum within the given time limit.). Column `time' reports the average computing time, when it is {not negligible}. The last row of the table reports the total number of optimal values found and the average gap and time, for each lower bound.
We do not report  on $LB_2^{II}$ since it does not improves upon $LB_2^I$, and uses some more computing time. Table \ref{tab:LBsALL-classes} gives the same information grouped by instance class: in this case each row refers to 150 instances.

\begin{table}[ht]
\footnotesize
\tabcolsep 3pt
\renewcommand
\arraystretch{1.05}
\begin{tabular}{lrrrrrrrrrrrrrrrrrrrrrrrrrrrrr}
\toprule
&  & &    \multicolumn{4}{c}{$LB_0$} & &          \multicolumn{3}{c}{$LB_0(\tilde{t})$} & &         \multicolumn{3}{c}{$LB_1$} & &          \multicolumn{3}{c}{$LB_2^I$} & &        \multicolumn{2}{c}{$LB_3$}\\
\cline{4-7} \cline{9-11} \cline{13-15} \cline{17-19} \cline{21-22}
       &          &&                  &         avg       &         max       &         avg       &&                  &         avg       &         max       &&                  &         avg       &         max       &&                  &         avg       &         max       &&                  &         avg       \\
$|T|$  &$|N|$     &&   \# opt    &         gap       &         gap       &         time      &&        \# opt    &         gap       &         gap       &&        \# opt    &         gap       &         gap       &&        \# opt    &         gap       &         gap       &&        \# opt    &         time      \\
\cline{1-2}\cline{4-7} \cline{9-11} \cline{13-15} \cline{17-19} \cline{21-22}
10     &54.90    &&97   &0.03      &1         &0.08&&      92        &0.08      &1         &&12        &1.49      &4         &&22        &1.20      &4      &&100       &0.24 \\
20     &88.39    &&96   &0.04      &1         &0.15&&      88        &0.13      &2         &&3         &1.82      &4         &&12        &1.46      &4      &&100       &0.90 \\
30     &121.43   &&95   &0.05      &1         &0.21&&      83        &0.18      &2         &&2         &2.02      &4         &&8         &1.66      &4      &&100       &2.23 \\
40     &154.10   &&96   &0.04      &1         &0.29&&      83        &0.18      &2         &&2         &2.09      &4         &&6         &1.72      &4      &&100       &4.23 \\
50     &186.70   &&98   &0.02      &1         &0.36&&      83        &0.21      &3         &&1         &2.31      &5         &&3         &1.90      &4      &&100       &5.51 \\
60     &219.61   &&97   &0.03      &1         &0.45&&      82        &0.22      &3         &&1         &2.27      &5         &&3         &1.89      &4      &&100       &9.80 \\
70     &252.50   &&98   &0.02      &1         &0.52&&      78        &0.26      &3         &&1         &2.31      &5         &&3         &1.93      &4      &&100       &14.41\\
80     &285.93   &&98   &0.02      &1         &0.57&&      82        &0.20      &3         &&1         &2.40      &5         &&5         &1.95      &4      &&100       &18.07\\
90     &318.86   &&98   &0.02      &1         &0.62&&      79        &0.24      &3         &&1         &2.49      &5         &&2         &2.04      &4      &&100       &30.14\\
100    &351.78   &&99   &0.01      &1         &0.69&&      77        &0.27      &3         &&1         &2.53      &5         &&2         &2.10      &4      &&100       &44.61\\
110    &385.04   &&98   &0.02      &1         &0.76&&      79        &0.27      &3         &&1         &2.55      &5         &&1         &2.10      &4      &&100       &29.79\\
120    &417.73   &&99   &0.01      &1         &0.84&&      80        &0.27      &3         &&0         &2.61      &5         &&0         &2.15      &4      &&100       &38.23\\
130    &451.12   &&99   &0.01      &1         &0.92&&      78        &0.30      &3         &&0         &2.58      &5         &&0         &2.14      &4      &&100       &48.40\\
140    &483.61   &&99   &0.01      &1         &0.97&&      79        &0.29      &3         &&0         &2.59      &5         &&1         &2.13      &4      &&100       &65.85\\
150    &516.45   &&99   &0.01      &1         &1.03&&      79        &0.30      &3         &&0         &2.61      &5         &&2         &2.14      &4      &&100       &67.47\\

\cline{4-7} \cline{9-11} \cline{13-15} \cline{17-19} \cline{21-22}	
       &          &&{\bf1466} &{\bf0.02}     &         &{\bf0.56}     &&{\bf1222}     &{\bf0.23}     &         &&{\bf26}      &{\bf2.31}     &         &&{\bf70}      &{\bf1.90}     &     &&{\bf1500}         &{\bf25.33}\\
\bottomrule
\end{tabular}
\caption{Performance of the Lower Bounds\label{tab:LBsALL-timeSteps}}
\end{table}
\begin{table}[ht]
\footnotesize
\tabcolsep 3pt
\renewcommand
\arraystretch{1.05}
\begin{tabular}{lrrrrrrrrrrrrrrrrrrrrrrrrrrrrr}
\toprule
&  & &    \multicolumn{4}{c}{$LB_0$} & &          \multicolumn{3}{c}{$LB_0(\tilde{t})$} & &         \multicolumn{3}{c}{$LB_1$} & &          \multicolumn{3}{c}{$LB_2^I$} & &        \multicolumn{2}{c}{$LB_3$}\\
\cline{4-7} \cline{9-11} \cline{13-15} \cline{17-19} \cline{21-22}
       &          &&                  &         avg       &         max       &         avg       &&                  &         avg       &         max       &&                  &         avg       &         max       &&                  &         avg       &         max       &&                  &         avg       \\
class  &$|N|$     &&   \# opt    &         gap       &         gap       &         time      &&        \# opt    &         gap       &         gap       &&        \# opt    &         gap       &         gap       &&        \# opt    &         gap       &         gap       &&        \# opt    &         time      \\
\cline{1-2}\cline{4-7} \cline{9-11} \cline{13-15} \cline{17-19} \cline{21-22}
I	&	89.00	&	&	144	&	0.04	&	1	&	0.16	&	&	135	&	0.10	&	1	&	&	13	&	1.35	&	2	&	&	27	&	1.03	&	2	&	&	150	&	0.49	\\
II	&	147.53	&	&	149	&	0.01	&	1	&	0.31	&	&	112	&	0.31	&	2	&	&	8	&	1.84	&	4	&	&	15	&	1.55	&	3	&	&	150	&	3.55	\\
III	&	164.94	&	&	150	&	0.00	&	0	&	0.32	&	&	90	&	0.49	&	3	&	&	0	&	2.45	&	4	&	&	2	&	2.16	&	4	&	&	150	&	3.14	\\
IV	&	208.98	&	&	147	&	0.02	&	1	&	0.52	&	&	118	&	0.22	&	2	&	&	1	&	2.15	&	4	&	&	4	&	1.57	&	4	&	&	150	&	17.45	\\
V	&	241.21	&	&	148	&	0.01	&	1	&	0.69	&	&	142	&	0.05	&	1	&	&	2	&	2.45	&	4	&	&	12	&	1.83	&	4	&	&	150	&	58.50	\\
VI	&	531.69	&	&	147	&	0.02	&	1	&	0.70	&	&	124	&	0.17	&	1	&	&	0	&	2.93	&	5	&	&	0	&	2.55	&	4	&	&	150	&	42.04	\\
VII	&	241.50	&	&	133	&	0.11	&	1	&	0.68	&	&	115	&	0.24	&	2	&	&	1	&	2.32	&	4	&	&	4	&	2.03	&	4	&	&	150	&	56.61	\\
VIII	&	326.86	&	&	148	&	0.01	&	1	&	0.75	&	&	121	&	0.25	&	2	&	&	1	&	2.49	&	4	&	&	4	&	1.87	&	3	&	&	150	&	17.39	\\
IX	&	554.91	&	&	150	&	0.00	&	0	&	0.69	&	&	133	&	0.13	&	2	&	&	0	&	2.55	&	4	&	&	2	&	2.21	&	3	&	&	150	&	23.13	\\
X	&	352.15	&	&	150	&	0.00	&	0	&	0.83	&	&	132	&	0.30	&	3	&	&	0	&	2.59	&	5	&	&	0	&	2.21	&	4	&	&	150	&	30.96	\\

\cline{4-7} \cline{9-11} \cline{13-15} \cline{17-19} \cline{21-22}	
       &          &&{\bf1466} &{\bf0.02}     &         &{\bf0.56}     &&{\bf1222}     &{\bf0.23}     &         &&{\bf26}      &{\bf2.31}     &         &&{\bf70}      &{\bf1.90}     &     &&{\bf1500}         &{\bf25.33}\\
\bottomrule
\end{tabular}
\caption{Performance of the Lower Bounds: instances grouped by groups \label{tab:LBsALL-classes}}
\end{table}
Bound $LB_0$ is quite effective and fails to find  the optimum solution value in only 34 instances over 1500. Its performances improve with the problem size.
Bound $LB_0(\tilde{t})$ fails on 278 instances and its performances definitely worsens when the size of the instance increases. Bounds $LB_1$ and $LB_2^I$ are not competitive since they are able to provide very few optimal values. The best bound is $LB_3$ which is able to give the optimal value for all the 1500 instances. However, this result is obtained with a larger computational effort. Indeed, it requires to solve the continuous relaxation \LPE\ for which the computing time grows  up to 67.47 seconds, on average,  with some rare instances were the time exceeds 2000 seconds.
Lower bound $LB_0$ runs, on average, in at most one second, while for $LB_0(\tilde{t})$, $LB_1$ and $LB_2^I$  the computing time is negligible and is not reported.

Note that although $LB_3$ equals  the optimal solution value in all the 1500 instances we generated, this equality cannot always hold if ${\cal P}\neq {\cal NP}$. For the \BPP\ it is known that the continuous relaxation of the exponential-size formulation almost always  provides the optimal solution value. The \emph{Integer Round-Up Property} (IRUP) states that the value of the LP relaxation of the exponential-size formulation, rounded up to the closest integer, yields the optimal solution value.
However,  for \BPP\ the  IRUP property  does not hold (see, e.g., \citet{CDDIR14}). The \TBPP\ generalizes the \BPP, so it also cannot exhibit the IRUP property.
\subsection{Heuristic algorithms}\label{subsec:comp-ub}
We implemented three heuristic algorithms. The first one, called ``Greedy" in the following, is made by running the four versions of the first fit method described in Section \ref{sec:ub} (First-Fit-1 and First-Fit-2 applied to the original and lifted instance) and  returning the best of these solutions.

The ``H-Rolling" heuristic implements the rolling horizon method introduced in Section \ref{sec:ub}. The method uses a parameter $\Delta$ to define the rolling time window. We performed preliminary tests  with $\Delta = \{10,20,30,40\}$, on a subset of instances. {On the basis of} these experiments, we selected the value $\Delta=30$ for our complete computational tests. We also made some preliminary tuning on the time limit given to \texttt{Cplex} for the solution of each restricted \ILPC\ model and we finally set the time limit to 10 seconds for instances with $|T|<100$ and to 30 seconds when $|T| \geq 100$.

The last heuristic ``H-Diving" is the diving method described in Section \ref{subsec:exp-heur}. Initially, it  solves the continuous relaxation \LPE with the column generation method described in Section \ref{sub:solvingLPe}, which provides both $LB_3$ and the starting point of the algorithm. Then, the branch-decision-tree defined by  branching rule BR-1 is partially explored  using the token-rule given in Section \ref{subsec:exp-heur}. For the computation of the root node ($LB_3)$ we set the time limit to 3500 CPU seconds, while 100 seconds are allowed for the tree exploration. The same  3600 seconds time limit will be used  for the exact algorithms.
\begin{table}[ht]
\footnotesize
\tabcolsep 4pt
\renewcommand
\arraystretch{1.05}
\begin{tabular}{crrrrrrrrrrrrrrrrrrrrrrrrrrrrr}
\toprule
&  & &    \multicolumn{4}{c}{Greedy} & & \multicolumn{4}{c}{H-Rolling} & & \multicolumn{4}{c}{H-Diving}\\
\cline{4-6} \cline{8-12} \cline{14-18}
          &           &&           & avg       & max       &&           & avg       & max       & avg    & max       &&           & avg       & max       & avg    & max          \\
$|T|$     & $|N|$     && \# opt    & gap       & gap       && \# opt    & gap       & gap       & time   & time      && \# opt    & gap       & gap       & time   & time         \\
\cline{1-2}\cline{4-6} \cline{8-12} \cline{14-18}
10       &54.90   &&46       &0.57     &2       &&100      &0.00     &0        &1.87     &10.00  && 100      &0.00    & 0        &0.46  &     7.02\\
20       &88.39   &&36       &0.84     &2       &&92       &0.08     &1        &4.63     &10.00  && 100      &0.00    & 0        &1.93  &    15.05\\
30       &121.43  &&32       &1.04     &3       &&79       &0.28     &3        &5.53     &10.00  && 100      &0.00    & 0        &5.92  &    65.20\\
40       &154.10  &&30       &1.10     &3       &&70       &0.34     &2        &5.95     &10.32  && 99       &0.02    & 2        &12.38 &   130.38\\
50       &186.70  &&36       &1.03     &3       &&75       &0.32     &3        &5.81     &20.00  && 97       &0.06    & 3        &18.93 &   155.72\\
60       &219.61  &&34       &1.06     &3       &&75       &0.29     &2        &6.88     &20.01  && 91       &0.15    & 3        &32.25 &   202.22\\
70       &252.50  &&29       &1.12     &3       &&72       &0.33     &2        &7.32     &20.08  && 87       &0.26    & 3        &44.49 &   299.81\\
80       &285.93  &&30       &1.13     &3       &&73       &0.30     &2        &7.41     &20.43  && 76       &0.47    & 3        &59.06 &   332.32\\
90       &318.86  &&30       &1.10     &4       &&77       &0.24     &2        &8.70     &25.77  && 81       &0.42    & 4        &71.37 &  1454.90\\
100      &351.78  &&27       &1.10     &4       &&78       &0.25     &2        &17.44    &60.68  && 69       &0.61    & 4        &98.22 &  2615.81\\
110      &385.04  &&27       &1.14     &3       &&79       &0.26     &3        &19.77    &65.81  && 63       &0.69    & 3        &84.16 &   580.28\\
120      &417.73  &&31       &1.09     &4       &&78       &0.26     &2        &20.06    &82.45  && 62       &0.70    & 4        &91.98 &  1272.17\\
130      &451.12  &&28       &1.14     &4       &&72       &0.33     &3        &19.57    &80.23  && 59       &0.80    & 4        &102.97&  1375.88\\
140      &483.61  &&31       &1.09     &3       &&74       &0.31     &2        &19.74    &78.42  && 50       &0.90    & 3        &127.66&  1701.14\\
150      &516.45  &&29       &1.14     &3       &&76       &0.31     &3        &19.68    &79.62  && 48       &0.93    & 3        &132.00&   602.14\\
\cline{4-6} \cline{8-12} \cline{14-18}
     &        &&{\bf476}      &{\bf1.05}     &           &&{\bf1170}     &{\bf0.26}     &{\bf}         &{\bf11.36}    &{\bf     }  && {\bf1182}     &{\bf0.40}    & {\bf}        &{\bf58.92} &   {\bf}\\
\bottomrule
\end{tabular}
\caption{Performance of the heuristic algorithms: instances grouped by time steps \label{tab:UBsALL-timeSteps}}
\end{table}
Table \ref{tab:UBsALL-timeSteps} gives the results for the three methods. The columns report on the number of times the heuristic solution is equal to the optimal solution (\# opt), the average and maximum absolute gap with respect to the optimum value (avg gap, max gap), and the average and max computing time, when it is relevant.
Each row provides results on the 100 instances we generated for each $|T|$ value. The last row summarises the above results on all the 1500 instances.
Table \ref{tab:UBsALL-Classes} provides the same information grouped by classes (we remind that in this case we have 150 instances per row).
\begin{table}[ht]
\footnotesize
\tabcolsep 4pt
\renewcommand
\arraystretch{1.05}
\begin{tabular}{lrrrrrrrrrrrrrrrrrrrrrrrrrrrrr}
\toprule
&  & &    \multicolumn{4}{c}{Greedy} & & \multicolumn{4}{c}{H-Rolling} & & \multicolumn{4}{c}{H-Diving}\\
\cline{4-6} \cline{8-12} \cline{14-18}
          &           &&           & avg       & max       &&           & avg       & max       & avg    & max       &&           & avg       & max       & avg    & max          \\
class     & $|N|$     && \# opt    & gap       & gap       && \# opt    & gap       & gap       & time   & time      && \# opt    & gap       & gap       & time   & time         \\
\cline{1-2}\cline{4-6} \cline{8-12} \cline{14-18}
I          &89.00   &&115      &0.24     &2       &&147      &0.02     &1        &0.05     &1.20    &&150      &0.00     &0        &1.41     &42.12  \\
II         &147.53  &&78       &0.49     &2       &&126      &0.16     &1        &1.89     &30.04   &&147      &0.02     &1        &16.45    &135.91 \\
III        &164.94  &&55       &0.66     &2       &&118      &0.21     &1        &2.16     &10.13   &&145      &0.05     &2        &14.16    &139.91 \\
IV         &208.98  &&14       &1.56     &3       &&123      &0.20     &2        &12.26    &40.77   &&107      &0.59     &3        &62.22    &246.13 \\
V          &241.21  &&13       &1.85     &4       &&74       &0.65     &3        &20.94    &79.62   &&102      &0.77     &4        &111.04   &2615.81\\
VI         &531.69  &&28       &1.06     &3       &&128      &0.17     &2        &14.53    &61.98   &&83       &0.67     &3        &100.78   &595.11 \\
VII        &241.50  &&6        &1.92     &3       &&108      &0.37     &3        &19.30    &82.45   &&106      &0.71     &3        &104.86   &1701.14\\
VIII       &326.86  &&44       &0.93     &3       &&120      &0.23     &3        &13.24    &60.63   &&115      &0.38     &3        &50.24    &428.01 \\
IX         &554.91  &&101      &0.36     &2       &&138      &0.08     &1        &9.26     &69.45   &&128      &0.16     &2        &48.60    &472.91 \\
X          &352.15  &&22       &1.39     &4       &&88       &0.52     &3        &19.93    &80.50   &&99       &0.66     &4        &79.43    &425.49 \\
\cline{4-6} \cline{8-12} \cline{14-18}
     &        &&{\bf476}      &{\bf1.05}     &           &&{\bf1170}     &{\bf0.26}     &{\bf}         &{\bf11.36}    &{\bf     }  && {\bf1182}     &{\bf0.40}    & {\bf}        &{\bf58.92} &   {\bf}\\
\bottomrule
\end{tabular}
\caption{Performance of the heuristic algorithms: instances grouped by class \label{tab:UBsALL-Classes}}
\end{table}
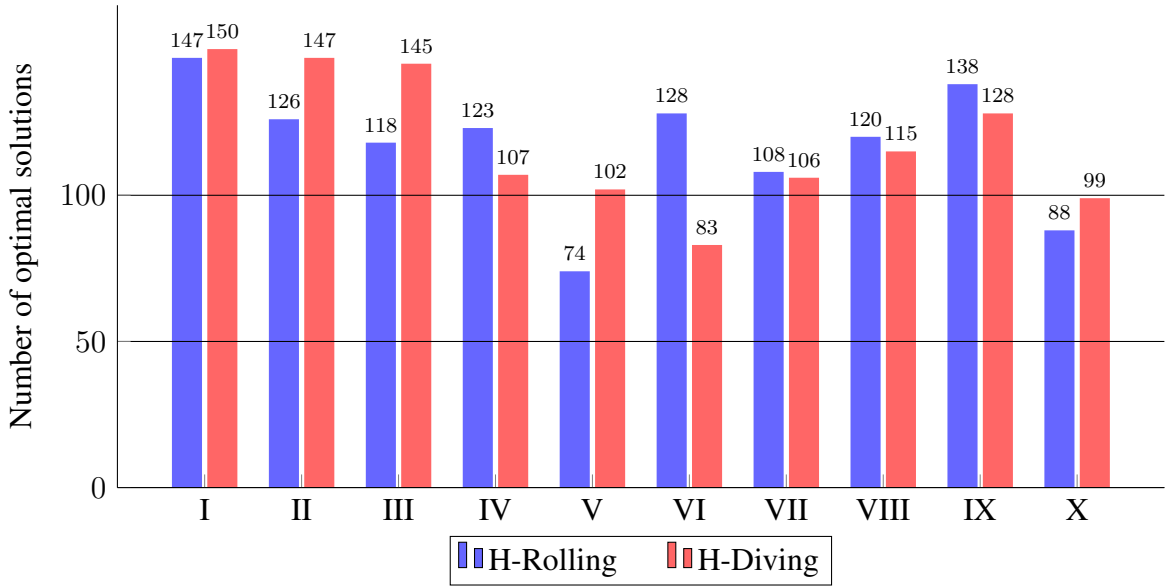
\begin{figure}[ht!]
\begin{center}
\pgfplotsset{width=7cm,compat=1.8}
\begin{tikzpicture}
  \centering
  \begin{axis}[
        ybar, axis on top,
        title={},
        height=8cm, width=15.5cm,
        bar width=0.4cm,
        ymajorgrids, tick align=inside,
        major grid style={draw=black},
        enlarge y limits={value=.1,upper},
        ymin=0, ymax=150,
        axis x line*=bottom,
        axis y line*=left,
        y axis line style={opacity=10},
        ytick={0, 50,  100}, tickwidth=5pt,
        enlarge x limits=true,
        legend style={
            at={(0.5,-0.1)},
            anchor=north,
            legend columns=-1,
            /tikz/every even column/.append style={column sep=0.6cm}
        },
        ylabel={Number of optimal solutions},
        symbolic x coords={
           I, II, III, IV, V, VI, VII, VIII, IX, X},
       xtick=data,
      nodes near coords={      \scriptsize
       \pgfmathprintnumber[precision=0]{\pgfplotspointmeta}
       }
    ]
    \addplot [draw=none, fill=blue!60] coordinates {
(I,     147)
(II,    126)
(III,   118)
(IV,    123)
(V,      74)
(VI,    128)
(VII,   108)
(VIII,  120)
(IX,    138)
(X,      88)
 };
    \addplot [draw=none, fill=red!60] coordinates {
(I,     150)
(II,    147)
(III,   145)
(IV,    107)
(V,     102)
(VI,    83 )
(VII,   106)
(VIII,  115)
(IX,    128)
(X,     99 )
 };
    \legend{H-Rolling, H-Diving}
  \end{axis}
  \end{tikzpicture}  \vspace*{-2ex}
\caption{Heuristic performances by instance classes}\label{fig:H-plot}
\end{center}
\end{figure}
{The Greedy is extremely fast and its computing time is not reported. It is able} to provide the optimal solution for about one third of the instances, while H-Rolling and H-Diving for about four over five solutions. H-Diving finds  twelve more optimal solutions than H-Rolling, but looking at each group of instances with the same $|T|$  one can see that no one of the two algorithms dominates the other.
If one looks at the single instances (not reported here) it is possible to find a sort of complementarity between the two methods: one solves instances not solved by the other and vice versa. The same considerations apply to the results grouped by class, shown in Table \ref{tab:UBsALL-Classes}. For example H-Diving improves upon H-Rolling for classes I, II, III, V and  X, while H-Rolling is the winner in the other classes. H-Rolling is able to find all the 100 optimal solutions only for $|T|=10$, then its performances slightly worsen, and for $|T|=(40,\dots,150)$ the number of optimal solutions found ranges between 70 and 80. H-Diving, instead, finds all the optimal solutions for $|T|\leq 30$, then its performances constantly decrease while $|T|$ increases. However, even when it finds only 48 optimal solutions ($|T|=150$) against the 76 found by H-Rolling, there are some instances where H-Diving beats H-Rolling (and vice versa).

From Table \ref{tab:UBsALL-Classes} and Figure \ref{fig:H-plot} one can see that class I is definitely the easiest for all algorithms, while classes II and III are also easy for H-Diving, although the computing times grows, on average, from 1.41 seconds to 16.45 seconds. Instances of classes IV and V are the most  difficult to solve for H-Rolling and  H-Divng, respectively.

\subsection{Exact algorithms}\label{subsec:comp-exact}
We started the computational analysis of the exact algorithms by comparing the results of the Vector Packing Solver by  \citet{Brandao2016}
with that of the polynomial-size formulation \eqref{eq:compact_obj}--\eqref{eq:compact_int_y}, solved by \texttt{Cplex} and with the exponential-size formulation solved by our branch-and-price algorithm introduced in Section \ref{sec:BP}. Table \ref{tab:small} reports the results for instances with $|T|=10, 15$. Each row refers to the 10 instances of a given {class}, and the columns show the results for the Vector Packing Solver {`VPSolver'}, for the polynomial-size model `\ILPC' and for the exponential-size formulation `\ILPE'.
For each algorithm, we provide the number of optimal solutions found (\# opt)  and the average computing time over the solved instances. The column labeled `$|N|$' reports the average number of items, while the column labeled `\#' reminds the number of instances tested in each class.
\begin{table}[ht]
\footnotesize
\tabcolsep 3pt
\renewcommand
\arraystretch{1.05}
\begin{tabular}{lrrrrrrrrrrrrrrrrrrrrrrrrrrrrrrrrrrrrrrrrrrrrrrrrrrr}
\toprule
                    &&&&\multicolumn{8}{c}{$|T|=10$}&&&\multicolumn{8}{c}{$|T|=15$}  \\
          \cline{5-12}\cline{15-22}
                    & &       & &       \multicolumn{2}{c}{{VPSolver}} &    & \multicolumn{2}{c}{\ILPC} & & \multicolumn{2}{c}{\ILPE}& & &        \multicolumn{2}{c}{{VPSolver}} &    & \multicolumn{2}{c}{\ILPC} & & \multicolumn{2}{c}{\ILPE}\\
\cline{5-6}\cline{8-9}\cline{11-12}\cline{15-16}\cline{18-19}\cline{21-22}

               class&         $|N|$     & \#    &         &         \# opt      &         time  &         &         \# opt      &         time  &         &         \# opt      &         time  &&&       \# opt      &         time  &         &         \# opt      &         time  &         &         \# opt      &         time\\
\cline{1-3}\cline{5-6}\cline{8-9}\cline{11-12}\cline{15-16}\cline{18-19}\cline{21-22}

          I         &         19.0      &         10        &         &        10         &      0.01         &         &                 10        & 0.02      &         &               10                  & 0.02      &  &        &             10        &    0.03      &          &                  10        &0.03      &             &         10 &         0.04         \\
          II        &         30.2      &         10        &         &        10         &      0.05         &         &                 10        & 0.06      &         &               10                  & 0.15      &  &        &             10        &    0.33      &          &                  10        &0.21      &             &         10 &         0.22         \\
          III       &         36.9      &         10        &         &        10         &      0.05         &         &                 10        & 0.14      &         &               10                  & 0.04      &  &        &             10        &    0.46      &          &                  10        &0.69      &             &         10 &         0.09         \\
          IV        &         46.1      &         10        &         &        10         &      0.33         &         &                 10        & 0.45      &         &               10                  & 0.71      &  &        &             10        &    3.66      &          &                  10        &1.98      &             &         10 &         1.10         \\
          V         &         53.8      &         10        &         &        10         &      0.48         &         &                 10        & 1.93      &         &               10                  & 0.47      &  &        &             10        &    9.70      &          &                  10        &6.20      &             &         10 &         1.59         \\
          VI        &         88.5      &         10        &         &        3          &      176.91       &         &                 10        & 3.91      &         &               10                  & 1.85      &  &        &             0         &    M.L.      &          &                  10        &48.79     &             &         10 &         6.09         \\
          VII       &         53.4      &         10        &         &        10         &      0.34         &         &                 10        & 0.47      &         &               10                  & 0.45      &  &        &             10        &    6.58      &          &                  10        &8.06      &             &         10 &         1.36         \\
          VIII      &         64.4      &         10        &         &        10         &      18.52        &         &                 10        & 0.91      &         &               10                  & 0.55      &  &        &             8         &    108.29    &          &                  9         &4.06      &             &         10 &         1.49         \\
          IX        &         87.6      &         10        &         &        4          &      193.75       &         &                 10        & 5.04      &         &               10                  & 1.13      &  &        &             0         &    M.L.      &          &                  10        &19.69     &             &         10 &         7.03         \\
          X         &         69.1      &         10        &         &        10         &      4.40         &         &                 10        & 125.17    &         &               10                  & 0.26      &  &        &             5         &   124.09     &          &                  10        &20.14     &             &         10 &         1.51         \\

\cline{5-6}\cline{8-9}\cline{11-12}\cline{15-16}\cline{18-19}\cline{21-22}

                    &                   &                   &         &    {\bf    87} & {\bf39.48}     &                  &  {\bf100}&       {\bf13.81}     & &              {\bf100}&            {\bf0.56}      &         & &    {\bf73} &         {\bf31.64}     &                & {\bf99}&         {\bf10.98}     &               &  {\bf100}&        {\bf2.05}      \\

\bottomrule
\end{tabular}
\caption{Performance comparison of exact algorithms on small instances (time limit 600 secs) \label{tab:small}}
\end{table}
A time limit of 600 seconds has been given to each algorithm. For $|T|=10$ ,the Vector Packing approach is not able to solve seven and six instances from classes VI and IX, respectively. For $|T|=15$, it is  able to solve only 73 over 100 instances. In this case, the unsolved instances are due to the excessive memory usage (entry `M.L.' in the table). The computing time drastically increases for some instances. {This is probably due to the fact that {VPSolver} is intended to solve a more general problem than the TBPP.} The \ILPC\ model is able to solve all instances with $|T|=15$ but one, using 31.64 CPU seconds on average. The \ILPE\ solves all instances in short times (on average 2.05 seconds) and is two order of magnitude faster than the {VPSolver}.
On the basis of these results, we decided that using {VPSolver} to solve the problem is not viable due to memory limitation ad high computing times, so we disregard this  approach for the next experiments.

We then performed computational experiments to assess the effectiveness of \ILPC\ and of our overall algorithm B\&P$^+$ to solve the \TBPP. We set a time limit of 3600 seconds for each algorithm and instance.
In Tables \ref{tab:ExactsALL-timeSteps} and \ref{tab:ExactsALL-classes} we report, for \ILPC\ and B\&P$^+$, the number of optimal solution and the average and maximum absolute gap and computing time. For B\&P$^+$ we additionally report the average root time, the average number of columns in the root LP, and the average number of nodes explored.
The average and maximum values are computed on the instances which are solved to a proven optimum. The two tables group the instances by $|T|$ and class, respectively.

From Table \ref{tab:ExactsALL-timeSteps} one can see that B\&P$^+$ dominates \ILPC\ for all $|T| < 150$, since it finds always more optimal solutions in a shorter average time. The \ILPC\ reaches the time limit for $|T|\geq 20$, while B\&P$^+$  for $|T|\geq 90$. For $|T|=150$ B\&P$^+$ finds two less optima than \ILPC\ in a comparable running time.
The analysis by instance classes in Table \ref{tab:ExactsALL-classes} confirms that B\&P$^+$ dominates \ILPC\ class by class, both by  number of optimal solutions found and  computing time. The first three classes appear to be easy for both algorithm, while class IV is also easy, but \ILPC\ fails on three instances. Class X remains the most difficult for both methods.

\begin{sidewaystable}[ht]
\footnotesize
\tabcolsep 5pt
\renewcommand
\arraystretch{1.05}
\begin{tabular}{lrrrrrrrrrrrrrrrrrrrrrrrrrrrrr}
\toprule
          &           &&    \multicolumn{6}{c}{\ILPC} & & \multicolumn{8}{c}{B\&P$^+$} \\
\cline{4-9} \cline{11-19}
          &           &&           &           &        &        &        &         &&           &           &           &        &              & \multicolumn{3}{c}{root}\\
\cline{16-19}
          &           &&           & avg       & max    & avg    & max    & avg     &&           & avg       & max       & avg    & max          & avg  & max   &      & avg     \\
$|T|$     & $|N|$     && \# opt    & gap       & gap    & time   & time   & nodes   && \# opt    & gap       & gap       & time   & time         & time  & time  & col  & nodes    \\
\cline{1-2}\cline{4-9} \cline{11-19}
10       &54.90   &&100      &0.00     &0.00     &13.81     &1234.59  &662.6      &&100      &0.00     &0        &1.95   &10.19    &0.24     &3.08     &103.1  &7.9\\
20       &88.39   &&99       &0.01     &1.00     &100.46    &t.l.     &2361.5     &&100      &0.00     &0        &5.24   &25.19    &0.90     &6.06     &264.7  &24.1\\                                                                                                               30       &121.43  &&96       &0.04     &1.00     &273.49    &t.l.     &2476.8     &&100      &0.00     &0        &9.33   &75.46    &2.23     &16.25    &524.9  &48.8\\                                                                                                           40       &154.10  &&92       &0.08     &1.40     &396.91    &t.l.     &2683.6     &&100      &0.00     &0        &15.74  &385.20   &4.23     &33.68    &844.7  &71.2\\                                                                                                           50       &186.70  &&93       &0.08     &1.39     &414.36    &t.l.     &2683.2     &&100      &0.00     &0        &27.89  &935.86   &5.51     &48.47    &1045.6 &85.9\\                                                                                                         60       &219.61  &&85       &0.17     &2.01     &596.10    &t.l.     &2289.5     &&100      &0.00     &0        &63.94  &2168.50  &9.80     &88.62    &1585.2 &114.9\\                                                                                                       70       &252.50  &&84       &0.18     &1.62     &694.06    &t.l.     &2648.7     &&100      &0.00     &0        &115.12 &2622.69  &14.41    &173.10   &2193.8 &148.1\\                                                                                                     80       &285.93  &&85       &0.20     &2.61     &635.66    &t.l.     &2627.2     &&100      &0.00     &0        &132.60 &2964.84  &18.07    &201.83   &2708.8 &174.8\\                                                                                                     90       &318.86  &&88       &0.15     &2.52     &578.75    &t.l.     &2955.5     &&99       &0.01     &1        &151.82 &t.l.     &30.14    &1295.94  &3178.4 &205.6\\                                                                                                     100      &351.78  &&86       &0.16     &2.00     &642.66    &t.l.     &2528.0     &&99       &0.01     &1        &168.85 &t.l.     &44.61    &2245.15  &3496.2 &207.8\\                                                                                                    110      &385.04  &&88       &0.14     &1.96     &609.58    &t.l.     &2754.6     &&98       &0.04     &3        &218.01 &t.l.     &29.79    &472.09   &4154.4 &229.4\\                                                                                                   120      &417.73  &&91       &0.12     &2.07     &523.89    &t.l.     &2353.4     &&99       &0.02     &2        &239.87 &t.l.     &38.23    &1164.56  &4588.4 &239.3\\                                                                                                  130      &451.12  &&92       &0.10     &2.02     &631.67    &t.l.     &2156.0     &&94       &0.09     &3        &433.91 &t.l.     &48.40    &1253.87  &4686.6 &238.3\\                                                                                                  140      &483.61  &&89       &0.14     &2.00     &634.70    &t.l.     &2687.4     &&92       &0.11     &2        &475.08 &t.l.     &65.85    &1555.42  &6764.3 &285.7\\                                                                                                            150      &516.45  &&92       &0.09     &2.34     &613.20    &t.l.     &2605.0     &&90       &0.15     &3        &560.40 &t.l.     &67.47    &502.21   &7519.6 &290.1\\                                                                                                              \cline{4-9} \cline{11-19}
         &        &&{\bf1360}&{\bf0.11}&{\bf}&          &         &{\bf2431.5}&&{\bf1471}&{\bf0.03}&         &{\bf174.75}&      &{\bf25.33} &     &{\bf2396.4}   &{\bf142.3}\\
\bottomrule
\end{tabular}
\caption{Performance of the exact algorithms, instances grouped by time steps} \label{tab:ExactsALL-timeSteps}
\end{sidewaystable}
\begin{sidewaystable}[ht]
\footnotesize
\tabcolsep 5pt
\renewcommand
\arraystretch{1.05}
\begin{tabular}{lrrrrrrrrrrrrrrrrrrrrrrrrrrrrr}
\toprule
          &           &&    \multicolumn{6}{c}{\ILPC} & & \multicolumn{8}{c}{B\&P$^+$} \\
\cline{4-9} \cline{11-19}
          &           &&           &           &        &        &        &         &&           &           &           &        &              & \multicolumn{3}{c}{root}\\
\cline{16-18}
          &           &&           & avg       & max    & avg    & max    & avg     &&           & avg       & max       & avg    & max          & avg  & max   &      & avg     \\
class     & $|N|$     && \# opt    & gap       & gap    & time   & time   & nodes   && \# opt    & gap       & gap       & time   & time         & time  & time  & col  & nodes    \\
\cline{1-2}\cline{4-9} \cline{11-19}
I         &      19.0 &&   150     & 0.00      &0.00    & 1.23   & 157.46 &  80.8   && 150       &0.00       &0          &0.55   &  30.69     &0.49     &6.81    &193.7    &11.2\\
II        &      30.2 &&   150     & 0.00      &0.00    & 3.90   & 82.68  &  107.3  && 150       &0.00       &0          &9.35   &  222.46    &3.55     &34.21   &809.9    &47.3\\
III       &      36.9 &&   150     & 0.00      &0.00    & 22.89  & 739.07 &  359.2  && 150       &0.00       &0          &8.10   &  279.36    &3.14     &37.83   &852.2    &65.7\\
IV        &      46.1 &&   147     & 0.02      &1.00    & 206.61 & t.l.   &  1444.9 && 150       &0.00       &0          &60.10  &  1135.04   &17.45    &141.83  &2222.9   &134.7\\
V         &      53.8 &&   128     & 0.19      &2.61    & 882.97 & t.l.   &  3449.9 && 146       &0.04       &2          &396.64 &  t.l.      &58.50    &2245.15 &2943.1   &164.1\\
VI        &      88.5 &&   142     & 0.06      &1.40    & 359.66 & t.l.   &  2066.4 && 146       &0.03       &1          &214.72 &  t.l.      &42.04    &473.19  &5238.6   &289.1\\
VII       &      53.4 &&   124     & 0.19      &2.34    & 874.38 & t.l.   &  2864.5 && 147       &0.04       &3          &249.37 &  t.l.      &56.61    &1555.42 &2499.9   &160.2\\
VIII      &      64.4 &&   132     & 0.15      &1.89    & 586.40 & t.l.   &  3256.9 && 144       &0.07       &3          &223.86 &  t.l.      &17.39    &326.13  &2770.3   &178.8\\
IX        &      87.6 &&   124     & 0.20      &2.00    & 845.39 & t.l.   &  4794.7 && 147       &0.02       &1          &82.28  &  t.l.      &23.13    &368.59  &2574.4   &132.7\\
X         &      69.1 &&   113     & 0.30      &2.07    & 1122.77& t.l.   &  5890.7 && 141       &0.09       &3          &501.58 &  t.l.      &30.96    &326.02  &3699.5   &229.7\\
\cline{1-2}\cline{4-9} \cline{11-19}
         &        &&{\bf1360}&{\bf0.11}&{\bf}&          &         &{\bf2431.5}&&{\bf1471}&{\bf0.03}&         &{\bf174.75}&      &{\bf25.33} &     &{\bf2396.4}   &{\bf142.3}\\
\bottomrule
\end{tabular}
\caption{Performance of the exact algorithms, instances grouped by classes} \label{tab:ExactsALL-classes}
\end{sidewaystable}
%

\begin{figure}[ht!]
\includegraphics[scale=1]{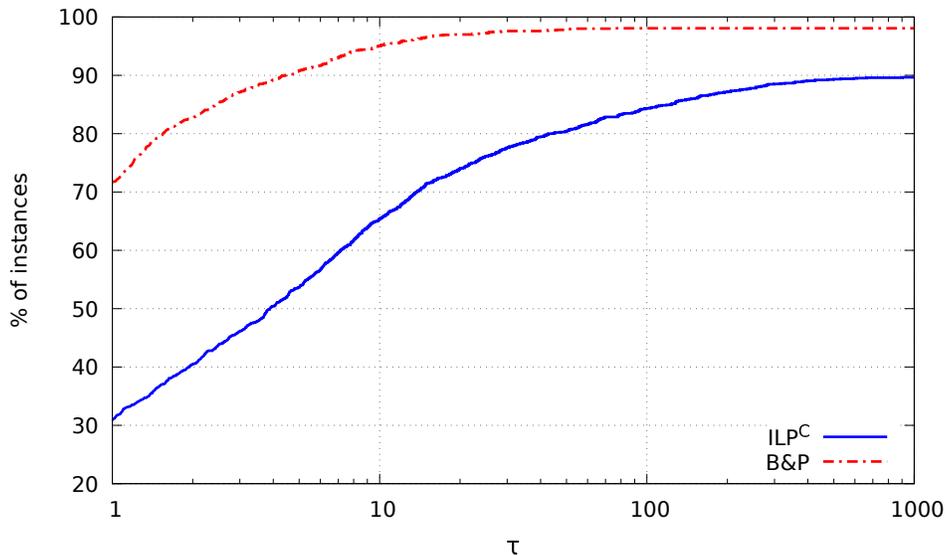}
\centering
\caption{Performance profile of the exact methods \label{fig:performance}}
\end{figure}

In order to give a graphical representation of the relative performance of the two exact algorithms, we report a performance profile in Figure \ref{fig:performance}. For each instance, we compute
a normalized time $\tau$ as the ratio of the computing time of the considered configuration over
the minimum computing time for solving the instance to optimality.
For each value of $\tau$ the vertical axis reports the percentage of the instances for which the corresponding
configuration spent at most $\tau$ times the computing time of the fastest configuration.
 The curves start from the percentage of instances in which the corresponding configuration is the fastest and
at the right end of the chart, we can read the percentage of instances solved by a specific algorithm.  The best performance are graphically represented by the curves in the upper part of  Figure \ref{fig:performance}. {B\&P$^+$ is the fastest one for approximatively 70\% of the instances and it is able to solve to problem optimality 98\% of them. On the other hand, \ILPC is able to solve only 90\% of the instances within the same time limit of 1 hour, and only 65\% and 85\% of the instances by allowing 10 and 100 times more time than the one required by B\&P$^+$, respectively. Figure \ref{fig:performance} graphically demonstrates that B\&P$^+$ compares favorably to \ILPC\ on the testbed of the 1500 considered instances. }

We complete our analysis by studying the contribution of each component of the overall algorithm B\&P$^+$ to the solution of the \TBPP, see
Tables \ref{tab:Contributions-timeSteps} and \ref{tab:Contributions-Classes}.
The first group of three columns (labeled H-Rolling+$LB_0$) refers to the first phase of the algorithms where we compute all the lower bounds, but $LB_3$ and we run the Greedy and the H-Rolling heuristic. In the  columns we report the number of instance solved (\# opt) and the average and maximum running times, computed with respect to the solved instances.
The second group of three columns (labeled H-Rolling+$LB_3$) refers to the instances that are not solved in the previous phase and shows the results obtained when lower bound $LB_3$
(i.e., the LP relaxation of the root node of the branch-decision-tree) is computed. Again we report the results evaluated only for the solved instances.
Note that in column labeled ``\#opt" we report the total number of instances solved in the first two phases of the algorithm, and in brackets the number of new optima.
The next group of three columns (labeled H-Diving) shows the results on the remaining unsolved instances, after the execution of the diving search.
The last group of five columns gives the results of the application of the B\&P$^+$ search to the instances that are not solved by the heuristics and lower bounds.
\begin{sidewaystable} 
\footnotesize
\tabcolsep 4pt
\renewcommand
\arraystretch{1.05}
\begin{tabular}{lrrrrrrcrrrcrrrrrrrrrrrrrrrrrr}
\toprule
          &           &&    \multicolumn{3}{c}{H-Rolling+$LB_0$}  &&\multicolumn{3}{c}{H-Rolling+$LB_3$}&& \multicolumn{3}{c}{H-Diving}&& \multicolumn{5}{c}{B\&P$^+$ search} \\
\cline{4-6} \cline{8-10} \cline{12-14}  \cline{16-20}
          &           &&              & avg    & max       &&           & avg    & max    &&           & avg    & max    &&           & avg       & max    & avg    & max  &          \\
$|T|$     & $|N|$     && \# opt       & time   & time      && \# opt    & time   & time   && \# opt    & time   & time   && \# opt    & gap       & gap    & time   & time &  \\
\cline{1-2} \cline{4-6} \cline{8-10} \cline{12-14}  \cline{16-20}
 10       &54.90   &&97       &1.80      &10.00 &&   100(3)   & 4.44     & 10.26  && -            &-         &-      &&  -       & -         & -        & -        & -\\ 20       &88.39   &&88       &3.90      &10.00 &&   92 (4)   & 12.35    & 13.72  &&100 (8)       &14.59     &24.39  &&  -       & -         & -        & -        & -\\
 30       &121.43  &&74       &4.08      &10.00 &&   79 (5)   & 14.37    & 22.10  &&100 (21)      &25.01     &72.27  &&  -       & -         & -        & -        & -\\
 40       &154.10  &&67       &4.49      &10.32 &&   70 (3)   & 18.51    & 31.74  && 99 (29)      &26.97     &93.46  &&100 (1)   & 0.00      &0.00      &385.20    &385.20  \\
 50       &186.70  &&73       &4.75      &20.00 &&   75 (2)   & 21.04    & 41.17  && 97 (22)      &34.33     &117.04 &&100 (3)   & 0.00      &0.00      &531.95    &935.86  \\
 60       &219.61  &&72       &6.10      &20.01 &&   75 (3)   & 20.55    & 58.11  && 95 (20)      &41.60     &106.51 &&100 (5)   & 0.00      &0.00      &999.49    &2168.50 \\
 70       &252.50  &&70       &6.45      &20.08 &&   72 (2)   & 56.69    & 73.83  && 94 (22)      &38.76     &98.22  &&100 (6)   & 0.00      &0.00      &1671.65   &2622.69 \\
 80       &285.93  &&72       &6.48      &20.43 &&   73 (1)   & 55.79    & 55.79  && 91 (18)      &44.97     &100.40 &&100 (9)   & 0.00      &0.00      &1318.28   &2964.84 \\
 90       &318.86  &&76       &8.02      &25.77 &&   77 (1)   & 76.80    & 76.80  && 92 (15)      &51.87     &116.66 &&99 (7)    & 0.13      &1.00      &1708.19   &t.l.    \\
 100      &351.78  &&78       &14.77     &60.68 &&   78 (-)   & -        & -      && 90 (12)      &66.79     &141.48 &&99 (9)    & 0.10      &1.00      &1486.10   &t.l.    \\
 110      &385.04  &&78       &17.82     &65.81 &&   79 (1)   & 93.17    & 93.17  && 88 (9)       &62.04     &126.55 &&98 (10)   & 0.33      &3.00      &1641.51   &t.l.    \\
 120      &417.73  &&78       &19.28     &82.45 &&   78 (-)   & -        & -      && 87 (9)       &57.14     &153.56 &&99 (12)   & 0.15      &2.00      &1684.90   &t.l.    \\
 130      &451.12  &&72       &17.70     &80.23 &&   72 (-)   & -        & -      && 83 (11)      &61.97     &123.74 &&94 (11)   & 0.53      &3.00      &2433.47   &t.l.    \\
 140      &483.61  &&74       &18.34     &78.42 &&   74 (-)   & -        & -      && 80 (6)       &52.06     &98.00  &&92 (12)   & 0.55      &2.00      &2288.29   &t.l.    \\
 150      &516.45  &&75       &16.67     &64.10 &&   76 (1)   &100.67    &100.67  && 82 (6)       &64.36     &114.64 &&90 (8)    & 0.83      &3.00      &3012.76   &t.l.    \\
 \cline{4-6} \cline{8-10} \cline{12-14}  \cline{16-20}
                   &                   &&{\bf1144}&{\bf9.93} &        &&{\bf1170 (26)}&{\bf28.22} &       &&{\bf1378 (208)}&{\bf41.85}&    && {\bf1471 (93)}&{\bf 0.35} &   &{\bf1970.21}\\ \bottomrule
\end{tabular}
\caption{Contribution of the various methods used in B\&P$^+$: instances grouped by time steps }\label{tab:Contributions-timeSteps}
\end{sidewaystable}

\begin{sidewaystable} \footnotesize
\tabcolsep 5pt
\renewcommand
\arraystretch{1.05}
\begin{tabular}{lrrrrrrcrrrrrrrrrrrrrrrrrrrrrr}
\toprule
          &           &&    \multicolumn{3}{c}{H-Rolling+$LB_0$}  &&\multicolumn{3}{c}{H-Rolling+$LB_3$}&& \multicolumn{3}{c}{H-Diving}&& \multicolumn{5}{c}{B\&P$^+$ search} \\
\cline{4-6} \cline{8-10} \cline{12-14}  \cline{16-20}
          &           &&              & avg    & max       &&           & avg    & max    &&           & avg    & max    &&           & avg       & max    & avg    & max  &          \\
$|T|$     & $|N|$     && \# opt       & time   & time      && \# opt    & time   & time   && \# opt    & time   & time   && \# opt    & gap       & gap    & time   & time &  \\
\cline{1-2} \cline{4-6} \cline{8-10} \cline{12-14}  \cline{16-20}
          I         &         19.0      &&141      &0.02      &1.20    &&14  (6)    &  1.19     &  2.43 && 150  (3 )     &16.44     &30.21  &&    -    &     -   &      -    &      -   & -   \\
          II        &         30.2      &&125      &2.09      &30.04   &&126 (1)    &  2.77     &  2.77 && 149  (23)     &37.58     &93.76  && 150 (1) & 0.00     & 0.00      &222.46    &222.46\\
          III       &         36.9      &&118      &2.26      &10.13   &&118 (-)    & -         & -     && 149  (31)     &20.02     &66.87  && 150 (1) &0.00     & 0.00      &279.36    &279.36\\
          IV        &         46.1      &&121      &12.32     &40.77   &&123 (2)    & 15.80     & 17.87 && 141  (18)     &46.96     &112.07 && 150 (9) &0.00     & 0.00      &729.38    &1135.04\\
          V         &         53.8      &&72       &18.41     &61.21   &&74  (2)    & 55.46     & 100.67&& 118  (44)     &51.63     &153.56 && 146 (28)&0.19     & 2.00      &1739.71   &t.l.   \\
          VI        &         88.5      &&125      &14.76     &61.98   &&128 (3)    & 57.38     & 76.80 && 136  (8 )     &53.52     &116.99 && 146 (10)&0.29     & 1.00      &2118.85   &t.l.   \\
          VII       &         53.4      &&98       &19.67     &82.45   &&108 (10)   & 38.43     & 93.17 && 130  (22)     &57.01     &131.99 && 147 (17)&0.30     & 3.00      &1687.26   &t.l.   \\
          VIII      &         64.4      &&118      &11.74     &41.10   &&120 (2)    & 12.37     & 13.59 && 138  (18)     &29.58     &98.22  && 144 (6) &0.83     & 3.00      &2627.64   &t.l.   \\
          IX        &         87.6      &&138      &8.12      &69.45   &&138 (-)    &  -        &   -   && 147  (9 )     &35.70     &72.27  && 147 (-) &1.00     & 1.00      &t.l.      &t.l.   \\
          X         &         69.1      &&88       &19.73     &80.50   &&88  (-)    &  -        &   -   && 120  (32)     &47.44     &140.60 && 141 (21)&0.47     & 3.00      &2396.26   &t.l.   \\
\cline{1-2} \cline{4-6} \cline{8-10} \cline{12-14}  \cline{16-20}
                    &                   &&{\bf1144}&{\bf9.93} &        &&{\bf1170 (26)}&{\bf28.22} &       &&{\bf1378 (208)}&{\bf41.85}&    && {\bf1471 (93)}&{\bf 0.35} &   &{\bf1970.21}\\
\bottomrule
\end{tabular}
\caption{Contribution of the various methods used in B\&P$^+$: instances grouped by classes} \label{tab:Contributions-Classes}
\end{sidewaystable}

\section{Conclusions}\label{sec:conclusions}
In this paper, {we studied the Temporal Bin Packing  Problem (TBPP), a challenging generalization of the classical Bin Packing Problem where each item consumes the bin capacity during a given time window. The goal is to determine the minimum number of bins to pack all the items while, at the same time, respecting the bin capacity at any instant of time. We have proposed and studied the first two mathematical formulations for the TBPP, the first one with a polynomial number of variables and constraints and the second one with an exponential number of variables. We have introduced several upper and lower bounds for the TBPP and we have designed an exact algorithm which combines them in an effective way. Our new branch-and-price algorithm, based on column generation, is able to solve to proven optimality instances with up to 500 items and 150 time steps,  in {reasonable computing times}.

Several are the potential future lines of research. In the recent literature, effective pseudo-polynomial size formulations have been proposed for the Bin Packing Problem. It would be interesting to study if these formulations, especially the ones based on the arc flow mechanism (see, e.g., \citet{V02} and \cite{DI18}), could be effectively used to tackle the TBPP as well. Finally, it would be interesting to introduce additional real-world features to the TBPP; for instance, precedence constraints between the items (see, e.g., \citet{DDI12}) or item class set-up costs (see, e.g., \citet{FMT18}).}
\section{Acknowledgments}
{Research supported by MIUR-Italy (Grant PRIN 2015, Nonlinear and Combinatorial Aspects of Complex Networks) and UNIMORE (Grant FAR 2018, Analysis and optimization of healthcare and pharmaceutical logistic processes).}
\clearpage
\section*{References}
\bibliographystyle{abbrvnat}
\bibliography{biblioTBPP}

\end{document}